\documentclass[11pt, dvips, twoside]{amsart}
\usepackage{amssymb,latexsym,bbm,graphicx,epsfig,epic,eepic,oldgerm}
\usepackage{eucal, amscd}
\usepackage{psfrag}
\usepackage{a4wide}
\usepackage{color}

\theoremstyle{plain}
\newtheorem{theorem}{Theorem}[section]
\newtheorem{lemma}[theorem]{Lemma}
\newtheorem{proposition}[theorem]{Proposition}
\newtheorem{corollary}[theorem]{Corollary}

\newtheorem{remarks}[theorem]{Remarks}

\newtheorem{question}[theorem]{Question}

\newcommand{\proofend}{\hspace*{\fill} $\Box$\\}

\def\s{\smallskip}
\def\m{\medskip}
\def\eps{\varepsilon}

\def\cat{\operatorname{cat}}

\def\cl{\operatorname{cl}}
\def\B{\operatorname {B}}
\def\C{\operatorname {C}}

\def\dim{\operatorname {dim}}

\def\Int{\operatorname {Int}\:\!}
\def\dist{\operatorname {dist}\:\!}

\def\diameter{\operatorname {diam}\:\!}

\def\st{\operatorname {st}}

\def\ot{\operatorname {ot}}
\def\rk{\operatorname {rk}}
\def\div{\operatorname {div}}

\def\size{\operatorname {size}}
\def\jet{\operatorname {jet}}
\def\rot{\operatorname {rot}}

\def\bb{\mathbf{b}}
\def\ccc{\mathbf{c}}
\def\xx{\mathbf{x}}
\def\yy{\mathbf{y}}
\def\00{\mathbf{0}}

\def\rr{\mathbf{r}}
\def\ss{\mathbf{s}}
\def\pphi{\mathbf{\phi}}

\def\ga{\alpha}
\def\gb{\beta}
\def\gg{\gamma}
\def\gd{\delta}
\def\eps{\varepsilon}

\def\gf{\varphi}

\def\gl{\lambda}
\def\go{\omega}
\def\gs{\sigma}

\def\cb{{\mathcal B}}
\def\cc{{\mathcal C}}
\def\cd{{\mathcal D}}

\def\cf{{\mathcal F}}

\def\ck{{\mathcal K}}
\def\cll{{\mathcal L}}
\def\cn{{\mathcal N}}

\def\cs{{\mathcal S}}

\def\cu{{\mathcal U}}
\def\cv{{\mathcal V}}

\def\ed{{\mathfrak D}}

\def\PP{\mathbbm{P}}

\def\RR{\mathbbm{R}}

\def\ZZ{\mathbbm{Z}}
\def\RP{{\RR\PP}} 
 
\def\11{\mathbbm{1}}

\def\pp{\partial}

\def\ra{\rightarrow}
\def\ha{\hookrightarrow}

\def\ni{\noindent}
\def\b{\bigskip}
\def\m{\medskip}

\def\im{\mbox{Im}\,}

\def\proof{\noindent {\it Proof. \;}}
\newcommand\proofof[1]{\noindent {\it Proof of #1. }}

\begin{document}

\numberwithin{equation}{section}



\title[]{Minimal atlases of closed contact manifolds}

\date{\today}

\author{Yuri Chekanov}
\address{(Yu.~Chekanov) 
Moscow Center for Continuous Mathematical Education,
B.~Vlasievsky per.~11, Moscow 119002, Russia}
\email{chekanov@mccme.ru}

\author{Otto van Koert}
\address{(O.~van Koert)
Department of Mathematics, Hokkaido University, Kita~10 Nishi~8, 
Sapporo 060-0810, Japan}
\email{okoert@math.sci.hokudai.ac.jp}

\author{Felix Schlenk}
\address{(F.~Schlenk) 
D\'epartement de Math\'ematiques, Universit\'e Libre de Bruxelles, CP~218, Boulevard du Triomphe, 1050 Bruxelles, Belgium} 
\email{fschlenk@ulb.ac.be}

\begin{abstract}  
We study the minimal number $\C (M,\xi)$ of contact charts that one needs to cover 
a closed connected contact manifold $(M,\xi)$.
Our basic result is $\C (M,\xi) \le \dim M + 1$.
We also compute $\C (M,\xi)$ for all closed connected contact $3$-manifolds:
\begin{eqnarray*}
\C (M,\xi) \,=\,
\left\{
 \begin{array}{ll} 
   2  & \text{if $M = S^3$ and $\xi$ is tight},  \\
   3  & \text{if $M = S^3$ and $\xi$ is overtwisted, or if $M = {\#}_k \left(S^2 \times S^1\right)$}, \\
   4  & \text{otherwise}.
 \end{array}
\right.
\end{eqnarray*}
We show that on every sphere $S^{2n+1}$ there exists 
a contact structure with \text{$\C (S^{2n+1}, \xi) \ge 3$}.
\end{abstract}

\maketitle

\begin{center}
\it
Dedicated to Yasha Eliashberg on the occasion of his 
sixtieth birthday
\end{center}

\section{Introduction and main results}  \label{s:intro}

\ni
A contact manifold is a pair $(M, \xi)$ where $M$ is a
smooth manifold of dimension $2n+1$ 
and $\xi \subset TM$ is a maximally non-integrable
field of hyperplanes. Such a field can be always 
written (at least locally) as the kernel of a  $1$-form~$\ga$. 
The maximal non-integrability condition then
has the form $\ga \wedge (d \ga)^n \neq 0$. 
The field $\xi$ is called a {\it contact structure}\/ on~$M$.
We refer to \cite{Gei-handbook,Gei-book,MS-98} 
for basic facts about contact manifolds.
The simplest contact manifold is $\RR^{2n+1}$ equipped with its
standard contact structure 
\[
\xi_{\st} \,=\, \ker \ga_{\st},  \quad \text{where } \ga_{\st} = dz + \sum_{i=1}^n x_i \,dy_i .
\]
A basic fact about contact manifolds is Darboux's Theorem which
states that locally every contact manifold $(M^{2n+1}, \xi)$ is contactomorphic 
to $(\RR^{2n+1}, \xi_{\st})$.
More precisely, for each point $p \in M$ there exists a chart 
$\phi \colon U \ra M$
from an open set $U \subset \RR^{2n+1}$
to $M$ such that $\phi^* \xi = \xi_{\st} |U$.
We call such a chart $\left( U, \phi \right)$ a {\it Darboux chart}.

\b
\ni
{\bf Definition.}
A {\it contact chart} for $(M,\xi)$ is a Darboux chart $\phi \colon U \to M$ 
such that $(U,\xi_{\st})$ is contactomorphic to $(\RR^{2n+1}, \xi_{\st})$.
The image $\phi (U)$ of a contact chart is a {\it contact ball}.

\b
\ni
As we shall see in Section~\ref{s:charts},
many subsets of $\RR^{2n+1}$ are contactomorphic to $(\RR^{2n+1}, \xi_{\st})$.
In particular, Eliashberg's classification of contact structures on 
$\RR^3$ in~\cite{El-92} implies that every subset of $\RR^3$ diffeomorphic to 
$\RR^3$ is contactomorphic to $(\RR^3, \xi_{\st})$.

\m
In the present paper we study the following

\m
\ni
{\bf Problem.}
{\it Given a closed contact manifold $(M, \xi)$, what is the minimal
number of contact charts that one needs to cover $(M, \xi)$?}

\m
\ni
In other words, we study the number $\C (M ,\xi)$ defined as 
\[
\C (M, \xi) \,=\, \min \left\{ k \mid M = \cu_1 \cup \dots \cup \cu_k  \right\}
\]
where each $\cu_i$ is a contact ball.
 
\m
An obvious lower bound for $\C (M, \xi)$ is the diffeomorphism
invariant 
\[
\B (M) \,=\, \min \left\{ k \mid M = B_1 \cup \dots \cup B_k \right\}
\]
where each $B_i$ is diffeomorphic to the standard open ball in $\RR^{2n+1}$. 

\m
Our basic result is

\m
\ni
{\bf Theorem 1.}
{\it Let $(M, \xi)$ be a closed connected contact manifold. 
Then 
$$
\B (M) \,\le\, \C (M, \xi) \,\le\, \dim M + 1.
$$
}

Apart from the trivial dimension~$1$, 
contact manifolds are best understood in dimension~$3$.
In this dimension, we compute $\C(M,\xi)$ for all closed contact manifolds.
If  $\xi$ is a contact structure on a $3$-manifold~$M$,
then the volume form $\ga \wedge d\ga$, 
where locally $\xi=\ker \ga$,
defines the {\it contact orientation}\/ $\nu_\xi$ on $M$. 
If $M$ already carries an orientation $\nu$, 
then the contact structure $\xi$ is called positive 
(with respect to $\nu$) if $\nu_\xi=\nu$, 
and negative otherwise.
Every closed oriented 
$3$-manifold admits a positive contact structure
in each homotopy class of tangent 
$2$-plane fields~\cite{El-89,Gei-handbook,Ma}.

Contact structures on $3$-manifolds fall into two classes, 
tight and overtwisted ones. This important dichotomy was introduced 
by Eliashberg~\cite{El-89}. The definitions go as follows.
A closed embedded $2$-disc $D$ in a contact $3$-manifold $(M,\xi)$
is called {\it overtwisted}\/ if $TD |_{\pp D} = \xi |_{\pp D}$.
A contact $3$-manifold $M$ is called {\it overtwisted}\/ if it contains 
an overtwisted disc, and {\it tight}\/ otherwise.
In this terminology, Bennequin's theorem~\cite{Ben} 
is equivalent to the tightness of the standard contact 
structure on the $3$-sphere~$S^3$. 
 
Overtwisted contact structures are more flexible,
their classification reduces to homotopy theoretical problems.
More precisely, it was proved by Eliashberg that
on an oriented closed \text{$3$-manifold},
every homotopy class of tangent $2$-plane fields contains 
a positive overtwisted contact structure, 
which is unique up to isotopy~\cite{El-89}.
On the other hand, tight contact structures are more rigid.
The first classification result here is again due to
Eliashberg, who showed that all tight contact structures on $S^3$ are 
contactomorphic. Since then, the classification of tight contact structures 
was achieved for many $3$-manifolds, 
see \cite{BC,Co-97,Co-99,Co-01,Co-02,CGH,CGH-preprint,El-92,Et-00,EG-99,EH-01,EH-02,Gh,GLS,GS,Gi-99,Gi-00,Gi-01,Go-98,Ho-00.I,Ho-00.II,Ho-02,Ho-04,HKM-02,HKM-03,HKM-04,Ka,KM,LM97,LM04,LS04,LS04b,Ma-Li}.

\m
\ni
{\bf Theorem 2.}
{\it
Let $(M,\xi)$ be a closed connected contact $3$-manifold $(M,\xi)$. Then
\begin{eqnarray*}
\C (M,\xi) \,=\,
\left\{
 \begin{array}{ll} 
   2  & \text{if $M = S^3$ and $\xi$ is tight},  \\
   3  & \text{if $M = S^3$ and $\xi$ is overtwisted, or if $M = {\#}_k \left(S^2 \times S^1\right)$}, \\
   4  & \text{otherwise}.
 \end{array}
\right.
\end{eqnarray*}
}

In other words, $\B (M) = \C (M,\xi)$ for all closed connected contact 
$3$-manifolds except for the 
countably infinite sequence $\xi_j$ of overtwisted contact structures 
on $S^3$, for which 
\text{$2 = \B (S^3) < \C (S^3, \xi_j) =3$}.
This shows that the contact invariant $\C(M,\xi)$ can be 
bigger than the smooth invariant $\B(M)$.
Theorem~2 solves in dimension~$3$ Problem~9.5 posed by Lutz 
in~\cite{Lu-88}:
$\C (S^3, \xi) =2$ if and only if $\xi$ is tight.

\b
The paper is organized as follows:
Section~\ref{s:top} provides methods for 
computing or estimating the lower bound $\B (M)$ 
of $\C (M,\xi)$.
In the rather technical Section~\ref{s:charts} we show that many subsets of 
$\RR^{2n+1}$ are contactomorphic to $(\RR^{2n+1}, \xi_{\st})$.
In Section~\ref{s:proof1} we prove Theorem~1, 
and in Sections~\ref{s:proof2}, \ref{s:convex}~and~\ref{s:proofend} 
we prove Theorem~2.
Section~\ref{s:ex} gives a few results on $\C(M,\xi)$ for contact 
manifolds of dimension $\ge 5$.

\b
\ni
{\bf Acknowledgments.}
We warmly thank 
Hansj\"org Geiges,
Emmanuel Giroux and
Klaus Niederkr\"uger 
for valuable discussions.
We in particular are grateful to 
Hansj\"org for showing us his new book~\cite{Gei-book} 
prior to publication,
and to
Klaus for contributing Proposition~\ref{p:PTN}.
Much of this work was done 
during 
the first authors stay at ULB in Spring 2007,
the second authors PostDoc at ULB in 2006 and 2007,
and 
the first and third authors stay at FIM of ETH Z\"urich in Summer 2008.
We wish to thank both institutions for their kind hospitality.

\section{Bounds for $\B(M)$}  \label{s:top}

\ni
Let $M$ be a smooth closed connected manifold of dimension $d$. 
Recall that $\B(M)$ denotes the minimal number of smooth balls 
covering $\B(M)$.
In view of Theorem~1 we are interested in 
lower bounds for $\B(M)$.

The {\it Lusternik--Schnirelmann category}\/ of $M$
is defined as
\[
\cat (M) \,=\, \min \{ k \mid M = A_1 \cup \ldots \cup A_k\}, 
\]
where each $A_i$ is open and contractible in $M$, 
\cite{CLOT,LuS}.
Clearly,
\[
\cat (M) \,\le\, \B (M) .
\]
There are examples with $\cat (M) < \B (M)$, 
see~\cite[Prop.~13]{LSV} and \cite[Prop.~3.6]{CLOT}.
In the case $d \le 3$, one always has $\cat (M) = \B (M)$,
cf.~Section~\ref{s:proof2} below.
In the case $d \ge 4$, sufficient conditions for 
$\cat (M) = \B (M)$ were found by Singhof~\cite{Si}.
It holds that
$\cat (M) = \cat (M')$ whenever $M$ and $M'$ are homotopy equivalent. 
The Lusternik--Schnirelmann
category is very different from the usual homotopical invariants
in algebraic topology and hence often difficult to compute.
Nevertheless, $\cat (M)$ can be estimated from below in 
cohomological terms as follows. 

The {\it cup-length}\/ of $M$ is defined as
\[
\cl (M) \,=\, \sup \{ k \mid u_1 \cdots u_k \neq 0, u_i \in
\tilde{H}^*(M)\} ,
\]
where $\tilde{H}^*(M)$ is the reduced cohomology of~$M$.
Then
\begin{equation*}  
\cl (M) +1 \le \cat (M),
\end{equation*}
see \cite{FE}. 
Given two closed connected manifolds~$M$ and $M'$, 
the LS-category of their 
product satisfies the following inequalities:
\begin{equation*}  
\max \{\cat (M), \cat (M') \} \le \cat ( M \times M') \le
 \cat (M) + \cat (M') -1.
\end{equation*}
Proofs of the above statements and more information on 
LS-category can be found in \cite{{CLOT}, {J1}, {J2}}.

Summarising, 
we have
\begin{equation} \label{e:EE}
 \cl (M) +1 \,\le\, \cat (M) \,\le\, \B (M)
\end{equation}
 for every closed connected manifold~$M$.
The upper bound
\begin{equation}  \label{e:Brough}
\B (M) \le d+1
\end{equation}
was proved in~\cite{L,OS},
and our proof of Theorem~1 will yield another proof.
Recall that $M$ is said to be $p$-connected 
if its homotopy groups 
$\pi_i (X)$ vanish for $1 \leq i \le p$. 
The following estimate considerably improves the 
estimate~\eqref{e:Brough}. 
\begin{proposition}  \label{p:Bupper}
Let $M$ be a closed connected smooth manifold of dimension $d \neq 4$.
If $M$ is $p$-connected then
\[
\B(M) \,\le\, \frac{d}{p+1} +1.
\]
\end{proposition}

\proof
For $d \ge 5$, this has been proved by Luft~\cite{L}, 
with the help of Zeeman's engulfing method~\cite{Z}.
The claim is obvious for $d \le 1$ and well-known for $d=2$.
For $d=3$, we can assume $p=1$ in view of~\eqref{e:Brough}, 
and then invoke the proof of the Poincar\'e conjecture.
\proofend

Proposition~\ref{p:Bupper} shows that Theorem~1 is far 
from sharp if $M$ is simply connected. 
The identity $\C (S^3,\xi_j) = 3$ from Theorem~2 shows that 
there is no analogue of Proposition~\ref{p:Bupper} for 
the contact covering number $\C(M,\xi)$.

A closed manifold $M$ with $\B (M) =2$ is homeomorphic to $S^d$, 
see~\cite{Br-60}.  
On the other hand, Proposition~\ref{p:Bupper} implies that 
a manifold $M$ homeomorphic to $S^d$ has $\B (M) =2$ 
provided that $d \neq 4$. 
We conclude that for a closed contact manifold,
\[
\B (M) =2 \,\Leftrightarrow\, \text{$M$ is homeomorphic to $S^d$}.
\]

\section{Contact charts}  \label{s:charts}

\ni
We shall often write $\RR^{2n+1}_{\st}$ for $\RR^{2n+1}$ 
endowed with its standard contact structure $\xi_{\st}$.
Two open subsets $U,V$ of $\RR^{2n+1}_{\st}$ are {\it contactomorphic}\/
if there exists a diffeomorphism $\gf \colon U \to V$ preserving $\xi_{\st}$,
that is, $d \gf (u) \xi_{\st} (u) = \xi_{\st}(\gf(u))$ for all $u \in U$.
Recall that a {\it contact chart}\/ for a contact manifold $(M,\xi)$ is 
a Darboux chart $\phi \colon U \to M$ 
where $U \subset \RR^{2n+1}_{\st}$ is contactomorphic to $\RR^{2n+1}_{\st}$.
In this section, we show that many subsets of $\RR^{2n+1}$ are contactomorphic 
to $(\RR^{2n+1}, \xi_{\st})$ and hence can be the domain $U$ of a 
contact chart.

A vector field $X$ on a contact 
manifold $(M,\xi)$
is called a {\it contact vector field}\/
if its local flow preserves~$\xi$.
An important example is the vector field
\[
V (\xx, \yy, z) \,=\, (\xx,\yy,2z)
\]
on $\RR^{2n+1}_{\st}$,
where $\xx,\yy \in \RR^n$, $z \in \RR$. 
The flow maps of $V$ are the contact dilations  
\[
\delta^t (\xx, \yy, z) \,=\, \left( e^t \xx, e^t\yy, e^{2t}z \right) .
\]
More generally, 
\begin{equation*}  
X (\xx, \yy, z) \,=\, (a\xx,b\yy,cz)
\end{equation*}
with $a,b,c \in \RR$ is a contact vector field iff $a+b=c$.

Let $(M,\xi=\ker\ga)$ be a contact manifold. 
There is a unique vector field $R_\ga$ on $M$ 
(the {\it Reeb vector field}\/) such that 
$i_{R_\ga}\ga = 1$ and $i_{R_\ga} d\ga = 0$.
For  $\RR^{2n+1}_{\st}$ one has $R_{\ga_{\st}}=(0,0,1)$.
The vector spaces of
contact vector fields on $(M,\xi)$  and
smooth functions on $M$ ({\it contact Hamiltonians}\/) 
are isomorphic via
\begin{equation}  \label{e:contactHam}
\left\{
 \begin{array}{lcl} 
    X &\mapsto& H_X = i_X \ga ;  \vspace{0.2em} \\
    H &\mapsto& X_H, \text{ with } 
    i_{X_H} \ga = H \text{ and } i_{X_H} d\ga = (i_{R_\ga}dH) \ga - dH .
 \end{array}
\right.
\end{equation}
For instance, the vector field $V(\xx,\yy,z) = (\xx,\yy,2z)$
corresponds to $H_{V}(\xx,\yy,z) = 2z +\xx\yy$
(where $\xx\yy=\sum_i x_i y_i$).
More generally, the vector field $X(\xx,\yy,z) = (a\xx,b\yy,cz)$ with $a+b=c$ corresponds to $H_X (\xx,\yy,z) = cz +b\:\!\xx\yy$.

A vector field $X$ on $\RR^{2n+1}$ is {\it complete}\/ if its flow $\gf_X^t$ exists for all times $t \in \RR$.
We say that a bounded domain $U \subset \RR^{2n+1}_{\st}$ is 
{\it contact star-shaped}\/ if there exists a complete contact vector field 
$X$ on $\RR^{2n+1}_{\st}$ such that
each flow line of $X$ intersects the boundary $\pp U$ 
in exactly one point
and such that
$$
\bigcup_{t \ge 0}\gf_X^t (U) \,=\, \RR^{2n+1}.
$$
The vector fields 
$X(\xx,\yy,z) = (a\xx,b\yy,cz)$ with $a,b \ge 0$ and $a+b=c >0$
provide many contact star-shaped domains containing the origin.

\begin{proposition}  \label{p:contact.starshaped}
Every  contact star-shaped domain $U \subset \RR^{2n+1}_{\st}$ is 
contactomorphic to $\RR^{2n+1}_{\st}$.
\end{proposition}

\proof
We follow Section~2.1 in~\cite{EG}.
Let $X$ be a contact vector field making $U$ a contact star-shaped domain,
and denote $U_t=\gf_X^t(U)$.
First we prove the following

\begin{lemma}  \label{abcd}
Given $a,b,c,d$ such that $a<b$ and  $c<d$, there is a 
contactomorphism $\Phi_{a,b}^{c,d}$ of $\RR^{2n+1}_{\st}$ that coincides 
with $\gf_X^{c-a}$ on a neighbourhood of $\partial U_a$ and
with $\gf_X^{d-b}$ on a neighbourhood of $\partial U_b$
(in particular, it sends $U_a$ to $U_c$ and $U_b$ to $U_d$). 
\end{lemma}

\proof
We can assume that $a=c$ 
(the general case can be reduced to this one by defining 
$\Phi_{a,b}^{c,d}=\gf_X^{c-a}\circ\Phi_{a,b}^{a,d-c+a}$).
Pick $q_1,q_2$ such that $a<q_1<q_2<\min \{ b,d \}$. 
Let $F$ be a smooth function on $\RR^{2n+1}$ such that
$F(u)=0$ when $u \in U_{q_1}$ and   
$F(u)=H_X$ when  $u \notin U_{q_2}$.
Then the contact vector field $X_F$ vanishes on  $U_{q_1}$
and coincides with $X$ outside of~$U_{q_2}$.
Its time $d-b$ flow map has the required properties.
\proofend

Choose two strictly increasing sequences of positive numbers, 
$(s_n)$ and  $(r_n)$,
such that $s_n\to 1$ and $r_n\to +\infty$.
The map $\varphi$ defined by 
\begin{eqnarray*}
\varphi (u) \,=\,
\left\{
 \begin{array}{ll} 
   \gf_X^{r_1-s_1}(u)  & \text{ when $u\in U_{s_1}$}, \\[0.2em]
     \Phi_{s_{n},s_{n+1}}^{r_{n},r_{n+1}}(u) & 
\text{ when $u\in U_{s_{n+1}}\setminus U_{s_n}$}
 \end{array}\right.
\end{eqnarray*}
is a contactomorphism from $U=U_1$ to $\RR^{2n+1}_{\st}$.
\proofend 

By a {\it cuboid}\/ in $\RR^{2n+1}$ we mean a bounded domain defined by hyperplanes parallel to the coordinate hyperplanes.

\begin{lemma}  \label{l:cuboid}
Every cuboid $Q$ in $\RR^{2n+1}_{\st}$ is contact star-shaped via a contact vector field vanishing in exactly one point in $Q$.
\end{lemma}

Together with Proposition~\ref{p:contact.starshaped} we obtain

\begin{corollary}  \label{c:cuboid}
Every cuboid in $\RR^{2n+1}_{\st}$ is contactomorphic to $\RR^{2n+1}_{\st}$.
\end{corollary}

\ni
{\it Proof of Lemma~\ref{l:cuboid}.}
Let $(\xx^0,\yy^0,z^0) \in \RR^{2n+1}$ be the centre of our cuboid
$$
Q \,=\, \times_{i=1}^n \left] x_i-a_i,x_i+a_i \right[ 
\times_{i=1}^n \left] y_i-b_i,y_i+b_i \right[ \times
\left] z-c,z+c \right[. 
$$
The affine contactomorphism
\begin{equation}  \label{e:tau}
\tau (\xx,\yy,z) \,=\,
\bigl( \xx-\xx^0,\yy-\yy^0,z-z^0+\xx^0 (\yy-\yy^0) \bigr)
\end{equation}
maps $(\xx^0,\yy^0,z^0)$ to the origin, and the faces of $\tau(Q)$ lie
in the $2n+1$ pairs of hyperplanes
\begin{eqnarray}
 &&\left\{ (\xx,\yy,z) \in \RR^{2n+1} \mid x_i= \pm a_i \right\},  \label{eq:cube_x_plane} \vspace{0.3em} \\
 &&\left\{ (\xx,\yy,z) \in \RR^{2n+1} \mid y_i= \pm b_i \right\},  \label{eq:cube_y_plane} \vspace{0.3em} \\
 &&\left\{ (\xx,\yy,z) \in \RR^{2n+1} \mid z = \pm c + \xx^0 \yy \right\} .  \label{eq:cube_z_plane}
\end{eqnarray}
As we have seen above, the vector field 
$$
X_\eps (\xx,\yy,z) \,:=\, \bigl( \eps \xx, (1+\eps)\yy, (1+2\eps)z \bigr)
$$
is a contact vector field for all $\eps >0$.
It is clearly transverse to the faces $\tau(\overline Q) \cap \eqref{eq:cube_x_plane}$ and $\tau(\overline Q) \cap \eqref{eq:cube_y_plane}$, pointing outward.
For $\eps >0$ small enough, 
$X_\eps$ is also transverse to the two faces 
$\tau(\overline Q) \cap \eqref{eq:cube_z_plane}$, pointing outward.
Indeed, the derivative of the function $f(\xx,\yy,z) = z - \xx^0 \yy$ along $X_\eps$ is
\begin{eqnarray*}
\cll_{X_\eps} f \,=\, df (X_\eps) &=& (dz-\xx^0d\yy)X_\eps \\
&=&(1+2\eps)z-(1+\eps)\xx^0\yy .
\end{eqnarray*}
At $(\xx,\yy,z) \in \eqref{eq:cube_z_plane}$ we therefore have
\begin{eqnarray*}
\cll_{X_\eps} f (\xx,\yy,z) &=& (1+2\eps)(\pm c + \xx^0\yy) - (1+\eps)\xx^0\yy \\
&=&  (1+2\eps)(\pm c) + \eps \xx^0\yy .
\end{eqnarray*}
Since $\yy$ is bounded on $\tau (\overline Q)$, we can choose $\eps >0$ so small that
this expression is positive on the face 
$\tau (\overline Q) \cap \left\{ z = +c + \xx^0\yy \right\} 
= \tau (\overline Q) \cap \left\{ f = +c \right\}$ 
and negative on the face 
$\tau (\overline Q) \cap \left\{ z = -c + \xx^0\yy \right\} 
= \tau (\overline Q) \cap \left\{ f = -c \right\}$.
Then $X_\eps$ points outward along these two faces.
The vector field $X_\eps$ therefore makes $\tau(Q)$ a contact star-shaped domain,
and so $(\tau^{-1})_* X_\eps$ makes $Q$ a contact star-shaped domain.
Since $X_\eps$ vanishes only in the origin, $(\tau^{-1})_* X_\eps$ vanishes only in the centre
$(\xx^0,\yy^0,z^0)$ of $Q$.
\proofend

For a further class of subsets of $(\RR^{2n+1},\xi_{\st})$ that are contactomorphic 
to $(\RR^{2n+1},\xi_{\st})$ we refer to Proposition~\ref{p:balls}.
The following question is open.

\begin{question} \label{q}
Let $U \subset \RR^{2n+1}$ be a subset diffeomorphic to $\RR^{2n+1}$,
where $2n+1\ge 5$.
Is it true that $(U,\xi_{\st})$ is always
contactomorphic to $(\RR^{2n+1},\xi_{\st})$?
\end{question}
In dimension~3, since every tight contact structure on $\RR^3$ is 
contactomorphic to $\xi_{\st}$ by~\cite{El-92}, we have
\begin{proposition}  \label{p:chart.3}
Every  domain $U \subset \RR^{3}_{\st}$ 
diffeomorphic to $\RR^3$ is contactomorphic to $\RR^3_{\st}$.
\end{proposition}

Our choice of definition for a contact chart is due to the preference  
to work with atlases having only one type of chart.
One may instead consider atlases consisting of Darboux charts
with domains that are only diffeomorphic to~$\RR^{2n+1}$. 
The corresponding contact covering numbers $\widetilde \C (M,\xi)$
satisfy the inequality $\widetilde \C (M,\xi) \le \C (M,\xi)$.
In dimension~3,   Proposition~\ref{p:chart.3} implies
$\widetilde \C (M,\xi) = \C (M,\xi)$.
If the answer to Question~\ref{q} is ``yes'' (which is rather unlikely), 
then always $\widetilde \C (M,\xi) = \C (M,\xi)$. 
All our results 
(except possibly Proposition~\ref{p:overtwisted.spheres}) 
remain valid for the invariant $\widetilde \C (M,\xi)$, 
since always $\B (M) \le \widetilde \C (M,\xi) \le \C (M,\xi)$ and since
$\widetilde \C (M,\xi) = \C (M,\xi)$ in dimension~$3$.

\b
The following result will be very useful later on.
\begin{proposition}  \label{p:one}
Let $\cu_1, \dots, \cu_k$ be disjoint contact balls in a contact 
manifold $(M,\xi)$,
and let $\ck$ be a compact subset of $\cu_1 \cup \dots \cup \cu_k$.
Then there exists a single contact ball covering~$\ck$.
\end{proposition}

\proof
Let $\phi_j \colon \RR^{2n+1} \to \cu_j$ be 
the contact chart for the contact ball  $\cu_j \subset M$.
Denote $\ck_j = \ck \cap \cu_j$ and $K_j = \phi^{-1}_j (\ck_j)$.
Pick $R$ such that the ball $B_R(0)\subset\RR^{2n+1}$ of radius $R$
contains~$K_j$.
Pick a smooth compactly supported function $f \colon \RR^{2n+1} \to \RR$ 
such that $f|_{B_R(0)}=1$. 
The contact vector field $X_{{-}f H_{V}}$ coincides with the 
contracting vector field $-V$ on~$B_R(0)$. 
Hence for $t\ge 0$ its time $t$ flow map $\gf^t$ coincides 
with $\delta^{{-}t}$ on~$B_R(0)$.
For each  $\eps >0$, there exists  $T>0$  such that 
$\delta^{{-}T}(K_j)\subset B_\eps(0)$.
Then $\gf^{T}$
sends $K_j$ into~$B_{\eps}(0)$.
Since $\gf^T$ is compactly supported, the map  
$\psi_j$ defined by
\begin{eqnarray*}
\psi_j (u) \,=\,
\left\{
 \begin{array}{ll} 
   \bigl( \phi_j \circ \gf^T \circ \phi_j^{-1} \bigr) (u)  & 
   \text{when $u \in \cu_j$},  \\
   \, u  & \text{when $u \notin \cu_j$},
 \end{array}
\right.
\end{eqnarray*}
is a contactomorphism of $(M,\xi)$ with support in~$\cu_j$. 
Because the contact balls $\cu_j$ are 
disjoint, the contactomorphism
\[
\Psi \,=\, \psi_k \circ \dots \circ \psi_1.
\]
coincides with  $\psi_j$ on~$\cu_j$. Therefore, $\Psi$ maps $\ck$
into 
\[
\cv \,=\, \bigcup_{j=1}^k\phi_j(B_\eps(0)).
\]
Pick a contact chart  $\phi \colon \RR^{2n+1}_{\st} \to \cu \subset M$
and a vector field $Y$ on $M$ such that its time $1$ flow map
sends each of the points $\phi_j(0)$ into~$\cu$.
Let $\Phi_Y^t$ be the flow of~$Y$.
Using contact Hamiltonians one easily constructs
a contact vector field $X$ coinciding with $Y$ on each of 
the paths
$\Phi_Y^t(\phi_j(0))$, $t\in[0,1]$.
Its time $1$ flow map $\Phi$ is a contactomorphism sending
each $\phi_j(0)$ into~$\cu$.
Choosing $\eps$ small enough ensures that 
 $\Phi$ maps $\cv$ into~$\cu$. 
Then $\Phi \circ \Psi$ maps $\ck$ into~$\cu$ and hence the contact chart 
$(\Phi \circ \Psi)^{-1} \circ \phi \colon \RR^{2n+1} \to M$ covers~$\ck$.
\proofend

\begin{figure}[ht] 
 \begin{center}
  \psfrag{1}{$\Psi$}
  \psfrag{2}{$\Phi$}
  \leavevmode\epsfbox{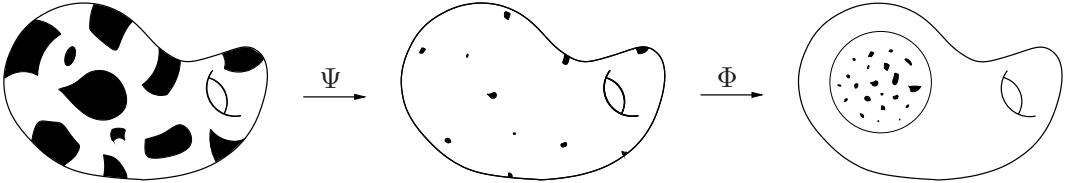}
 \end{center}
 \caption{The proof of Proposition~\ref{p:one}.} 
 \label{figure24}
\end{figure}
%
%

\section{Proof of Theorem~1}  \label{s:proof1}

\ni
We want to show that $\C (M,\xi) \le d+1$ for every closed connected contact 
manifold $(M,\xi)$ of dimension $d := 2n+1$.
An analogous result for symplectic manifolds is proved in~\cite{RS}.
We shall prove Theorem~1 by using the same idea as in~\cite{RS}.
The proof for contact covers is, however, much easier than the one for symplectic covers,
thanks to Proposition~\ref{p:one}.

\begin{figure}[ht] 
 \begin{center}
  \leavevmode\epsfbox{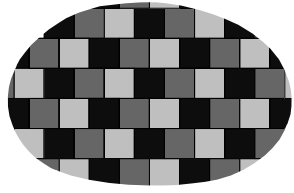}
 \end{center}
 \caption{A part of the dimension cover of $\RR^2$.} 
 \label{figure20b}
\end{figure}

\m
\ni
{\bf Idea of the proof.}
As one knows from looking at a brick wall,
the plane $\RR^2$ can be divided into squares
coloured with three colours in such a way that 
squares of the same colour do not touch, 
see Figure~\ref{figure20b}.

%
%

\ni
It is the starting-point of dimension theory, \cite{En},
that this observation extends to all dimensions:
$\RR^d$ can be divided into $d$-dimensional cubes
coloured with $d+1$ colours in such a way that cubes of the same colour do not touch,
see below for the construction. 
We shall show that this construction extends, in a way, to 
every closed contact manifold $(M,\xi)$ of dimension~$d$:
the manifold $M$ can be covered by $d+1$ sets $\cs^1, \dots, \cs^{d+1}$ 
such that each component of $\cs^j$ is contactomorphic to a cube in~$\RR^d$.
Corollary~\ref{c:cuboid} and Proposition~\ref{p:one} now imply that $\cs^j$ 
can be covered by a single contact ball, and so $M$ can be covered by $d+1$ contact balls.

\m
We now come to the actual proof of Theorem~1.

\b
\ni
{\bf 1. The cover by contact cubes}

\m
\ni
We start by constructing the standard dimension cover of $\RR^d$.
We do this by successively constructing certain covers of $\RR^i$, 
$i = 1, \dots, d$. 
Cover $\RR^1$ by intervals of the form $[k-1,k]$.
The cubes ``of colour $j$'' in this partition are the intervals
$$
\coprod_{k \in j+(d+1)\ZZ} [k-1,k] \,\subset\, \RR .
$$
For $2 \le i \le d$ the cover of $\RR^i = \{ (x_1,\dots,x_i) \}$ 
by cubes is obtained from the cover of $\RR^{i-1} = \{ (x_1,\dots,x_{i-1}) \}$ 
by the following procedure. 
The first layer $\{ 0\le x_i \le 1 \} \subset \RR^i$ is filled by $i$-dimensional cubes 
whose projection to $\RR^{i-1}$ form the partition of~$\RR^{i-1}$,
and the colour-$j$-cubes in the first layer are those that project to the 
colour-$j$-cubes in the partition of $\RR^{i-1}$.
We construct the cubes of the $k$-th layer $\{ k-1 \le x_i \le k\}$, $k \in \ZZ$,
by applying the translation by $\bigl( (1/d)(k-1) ,\ldots, (1/d) (k-1), k-1 \bigr)$ 
to the cubes of the first layer.
The colour-$j$-cubes on the $k$-th layer are obtained from 
the colour-$j$-cubes of the first layer by applying the translation 
$\bigl( (1+1/d)(k-1) ,\ldots, (1+1/d) (k-1), k-1 \bigr)$.

\begin{figure}[ht] 
 \begin{center}
  \psfrag{x1}{$x_1$}
  \psfrag{x2}{$x_2$}
  \leavevmode\epsfbox{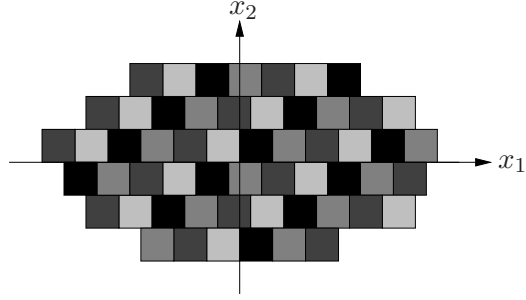}
 \end{center}
 \caption{A part of an $x_3$-layer of the dimension cover of $\RR^3$.} 
 \label{figure.dim.cover}
\end{figure}

Let $\ed^j(d)$ be the set in $\RR^d$ formed by the cubes of colour $j$.
The sets $\ed^1(d), \dots, \ed^{d+1}(d)$ form the standard dimension cover of $\RR^d$.
Figure~\ref{figure20b} shows a part of this cover for $d=2$,
and Figure~\ref{figure.dim.cover} shows a part of an $x_3$-layer of this cover 
for $d=3$.

%
%

\begin{figure}[ht] 
 \begin{center}
  \leavevmode\epsfbox{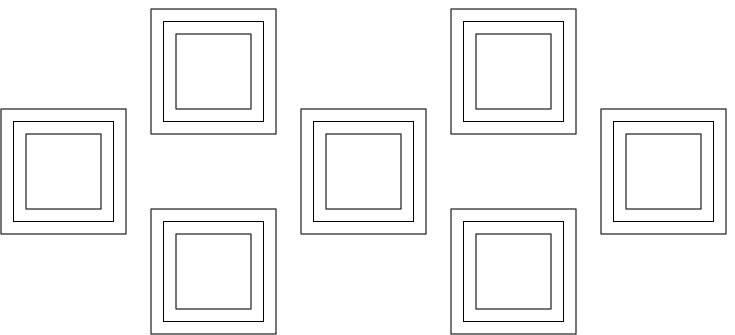}
 \end{center}
 \caption{Seven cubes in $s \ed^j(2)$ and their neighbourhoods $\cn_1(C)$ and $\cn_2(C)$.} 
 \label{figure.cn}
\end{figure}

For each $s>0$ we have a dimension cover of $\RR^d$ formed by the sets
$$
s \ed^j(d), \quad \; j = 1, \dots, d+1 ,
$$
which are the images of $\ed^j(d)$ under the homothety $v \mapsto s v$ of $\RR^d$.
Given a cube $C$ in $s \ed^j(d)$ we denote by $\cn_1(C)$ and $\cn_2(C)$
the closed cubes in $\RR^d$ with edges parallel to the axes 
that have the same centre as $C$ and have sizes (i.e.\ edge lengths) 
\begin{equation}  \label{e:size}
\size \left( \cn_1 (C) \right)  
\,=\, 
\left( 1+\frac{1}{4d} \right) \size (C) 
\,=\, 
\left( 1+\frac{1}{4d} \right) s
, \qquad
\size \left( \cn_2 (C) \right)  
\,=\, 
\left( 1+\frac{1}{2d} \right) s .
\end{equation}
In view of the construction of the standard dimension cover, 
the distance between a cube $C$ in $s \ed^j(d)$ and 
$s \ed^j(d) \setminus C$ is $\frac 1d s$, cf.~Figures~\ref{figure20b} and \ref{figure.dim.cover}.
Therefore, the neighbourhoods $\cn_2 (C)$ and $\cn_2 (C')$ of different cubes $C,C'$ in $s \ed^j(d)$ are disjoint,
see Figure~\ref{figure.cn}.

%
%

\m
Fix now $d=2n+1 = \dim M$, 
and let $B^d(1)$ be the open ball in $\RR^d$ of radius $1$ centred at the origin.
Since $M$ is compact, we can choose finitely many contact charts $\phi_k \colon \RR^d \to M$,
$k=1,\dots,\ell$, such that the $\ell$ open sets $\phi_k \left( B^d(1) \right)$
cover $M$.
Fix a colour $j$.
For $s>0$ we denote by $(s \ed^j(d))_1$ the set formed by the cubes in $s \ed^j(d)$ that intersect $B^d(1)$.
We are going to choose for each $k$ a size $s_k$ and a number $\gd_k$ in an appropriate way.
To this end, fix a distance on $M$ induced by a Riemannian metric on $M$.
First choose $s_\ell = 1$, 
and choose $\gd_\ell >0$ with 
$$
\gd_\ell \,\le\, \dist  \bigl( \phi_\ell(C), M \setminus \phi_\ell  (\cn_1(C)) \bigr),
\quad
\gd_\ell \,\le\, \dist  \bigl( \phi_\ell(\cn_1(C)), M \setminus \phi_\ell  (\cn_2(C)) \bigr)
$$
for each of the finitely many cubes $C$ in $(s_\ell \ed^j(d))_1$.
Assume by induction that we have chosen $s_\ell, s_{\ell-1}, \dots, s_{i+1} >0$ and 
$\gd_\ell \ge \gd_{\ell -1} \ge \dots \ge \gd_{i+1} >0$.
For $s \in \left] 0,1\right]$ and for each cube $C$ in $(s \ed^j(d))_1$ we have 
$$
\cn_2(C) \,\subset\, B^d \left( 1+(1+\tfrac{1}{2d})\sqrt d \right) .
$$
Since the differential of $\phi_i$ 
(with respect to the Euclidean metric on $\RR^d$ and the Riemannian metric on $M$)
is uniformly bounded on this ball,
we can choose $s_i \in \left] 0,1\right]$ such that for each cube $C$ in $(s_i \ed^j(d))_1$,
\begin{equation}  \label{e:ddk}
\diameter \bigl( \phi_i (\cn_2(C)) \bigr) \,<\, \gd_{i+1} .
\end{equation}
Now choose $\gd_i > 0$ so small that $\gd_i \le \gd_{i+1}$ and 
\begin{equation}  \label{e:dk}
\gd_i \,\le\, \dist  \bigl( \phi_i(C), M \setminus \phi_i  (\cn_1(C)) \bigr),
\quad
\gd_i \,\le\, \dist  \bigl( \phi_i(\cn_1(C)), M \setminus \phi_i (\cn_2(C)) \bigr)
\end{equation}
for each of the finitely many cubes $C$ in $(s_i \ed^j(d))_1$.

A {\it coordinate cube}\/ in $\RR^d$ is a closed cube 
with edges parallel to the axes.
A {\it contact cube}\/ in a contact manifold~$P$ is the image of a coordinate cube 
in $\RR^d$ under a contact chart $\RR^d \to P$.
\begin{proposition}  \label{p:new.cubes}
The set $\bigcup_{k=1}^\ell \phi_k \bigl( (s_k \ed^j(d))_1 \bigr)$ of cubes of colour $j$ 
can be covered by a finite disjoint union of contact cubes. 
\end{proposition}

\proof
We start with

\begin{lemma}  \label{l:onecube}
Let $C$ be a cube of $(s \ed^j(d))_1$, let
$\ck_1, \dots, \ck_a$ be disjoint contact cubes in $\Int \cn_2(C)$
that intersect $\cn_1(C)$ and are disjoint from $C$,
and let $U \subset \cn_2(C) \setminus C$ be an open neighbourhood of $\ck_1 \cup \dots \cup \ck_a$.
Then there exists a contactomorphism $\psi$ of $\RR^d_{\st}$ with support in $U$ 
such that $\psi (\cn_1(C)) \supset C \cup \ck_1 \cup \dots \cup \ck_a$.
\end{lemma}

\begin{figure}[ht] 
 \begin{center}
  \psfrag{C}{$C$}
  \psfrag{cn1}{$\cn_1(C)$}
  \psfrag{cn2}{$\cn_2(C)$}
  \psfrag{p}{$\psi$}
  \leavevmode\epsfbox{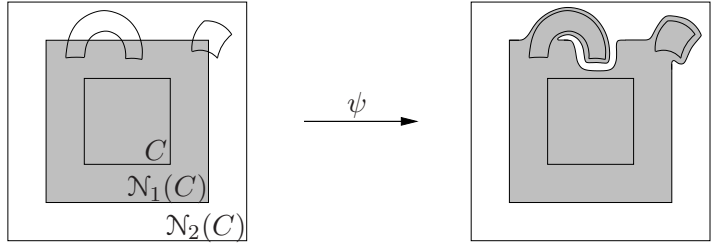}
 \end{center}
 \caption{The contactomorphism $\psi$.} 
 \label{figure.psi}
\end{figure}

\proof
Let $\Phi_\ga \colon \RR^d \to \im \Phi_\ga$ be a contact chart with $\Phi_\ga(K_\ga) = \ck_\ga$,
where $K_\ga$ is a coordinate cube in $\RR^d$.
After replacing $K_\ga$ by a slightly larger coordinate cube, if necessary,
we can assume that $\ck_\ga$ intersects $\Int \cn_1(C)$.
By Lemma~\ref{l:cuboid},
there exists a contact vector field on $\RR^d$ which vanishes at exactly one point in $\Int K_\ga$ and is transverse to $\pp K_\ga$.
Its image $X_\ga$ under $\Phi_\ga$ is a contact vector field on $\im \Phi_\ga$ which vanishes at exactly one point $p_\ga$ in $\Int \ck_\ga$ and is transverse to $\pp \ck_\ga$.
After applying a contactomorphism with support in $\ck_\ga$ we can assume that $p_\ga \in \Int \ck_\ga \cap \Int \cn_1(C)$.
Let $H_\ga$ be the contact Hamiltonian for $X_\ga$.
For each $\ga$ choose $\cu_\ga \subset U \cap \im \Phi_\ga$ such that $\cu_\ga \supset \ck_\ga$ 
and such that the sets $\cu_\ga$ are mutually disjoint,  
and choose a smooth function $f_\ga$ with support in $\cu_\ga$ such that $f_\ga |_{\ck_\ga} \equiv 1$.
The flow $\gf_\ga^t$ of the contact vector field $X_{f_\ga H_\ga}$ is then supported in $\cu_\ga$,
and for large enough $T_\ga$ we have $\gf_\ga^{T_\ga} \left( \Int \cn_1(C) \cup \cu_\ga \right) \supset \ck_\ga$,
cf.~Figure~\ref{figure.psi}.
Since $\cu_1, \dots, \cu_a$ are mutually disjoint and disjoint from $C$, the contactomorphism 
$\psi = \gf_a^{T_a} \circ \dots \circ \gf_1^{T_1}$ has the required properties.
\proofend

%
%

\ni
{\it Proof of Proposition~\ref{p:new.cubes}.}
We shall prove by induction that for each $i=1, \dots, \ell$,

\m
\ni
{\bf Claim~($i$).}
{\it The set $\bigcup_{k=1}^i \phi_k \bigl( (s_k \ed^j(d))_1 \bigr)$ can be covered by a finite disjoint union of contact cubes of diameter $< \gd_{i+1}$.}

\m
\ni
Claim~($\ell$) implies Proposition~\ref{p:new.cubes}.
Claim~(1) is obvious, since the components of $\phi_1 \bigl( (s_1 \ed^j(d))_1 \bigr)$ are contact cubes of diameter $< \gd_2$ according to~\eqref{e:ddk}.
Assume that Claim~($i-1$) holds true, that is, 
$\bigcup_{k=1}^{i-1} \phi_k \bigl( (s_k \ed^j(d))_1 \bigr)$ is covered by contact cubes
$\ck_1, \dots, \ck_a$ of diameter $< \gd_i$.
We want to cover $\ck_1, \dots, \ck_a$ and the contact cubes $\phi_i(C_1) = \cc_1, \dots, \phi_i(C_b) = \cc_b$ 
from $\phi_i \bigl( (s_i \ed^j(d))_1 \bigr)$ by disjoint contact cubes of diameter~$< \gd_{i+1}$. 
Write
$$
\cn_1(\cc_\gb) = \phi_i \left( \cn_1(C_\gb) \right)
\qquad \text{and} \qquad
\cn_2(\cc_\gb) = \phi_i \left( \cn_2(C_\gb) \right) .
$$
Since $\cn_2(C_\gb) \cap \cn_2(C_{\gb'}) = \emptyset$ for $\gb \neq \gb'$, we have 
\begin{equation}  \label{e:s'}
\cn_2(\cc_\gb) \cap \cn_2(\cc_{\gb'}) = \emptyset
\quad \text{for }\, \gb \neq \gb' . 
\end{equation}
Given $\gb \in \{ 1, \dots, b \}$, consider those contact cubes among 
$\ck_1, \dots, \ck_a$ that intersect $\cn_2(\cc_\gb)$.
By~\eqref{e:ddk} and \eqref{e:dk}, each of these contact cubes 
\begin{itemize}
\item[(i$_\gb$)]
either is contained in $\cn_1(\cc_\gb)$;
\item[(ii$_\gb$)]
or is contained in $\cn_2(\cc_\gb) \setminus \cc_\gb$;
\item[(iii$_\gb$)]
or is disjoint from $\cn_1(\cc_\gb)$.
\end{itemize}
We apply Lemma~\ref{l:onecube} to $C_\gb$, 
to the contact cubes $\phi_i^{-1}(\ck_\ga)$ with $\ck_\ga$ of type~(ii$_\gb$),  
and to a neighbourhood $U \subset \cn_2(C_\gb)$ which is disjoint from $C_\gb$ and from
$\phi_i^{-1}(\ck_\ga)$ for all $\ck_\ga$ of type~(i$_\gb$) or (iii$_\gb$).
We then obtain a contact cube $\ck_{\cc_\gb}$ which covers $\cc_\gb$ and the contact cubes of type~(i$_\gb$) and (ii$_\gb$).
Moreover, $\ck_{\cc_\gb}$ is contained in $\cn_2(\cc_\gb)$ and is disjoint from the cubes $\ck_\ga$ of type~(iii$_\gb$).
Since $\ck_{\cc_\gb} \subset \cn_2(\cc_\gb)$, the contact cubes $\ck_{\cc_\gb}$ are disjoint 
by~\eqref{e:s'},
and have diameter $< \gd_{i+1}$ by~\eqref{e:ddk}.
Let $\ck_{\ga_1}, \dots, \ck_{\ga_m}$ be those contact cubes among $\ck_1, \dots, \ck_a$
that are disjoint from $\bigcup_{\gb=1}^b \cn_1(\cc_\gb)$.
These contact cubes are disjoint from $\ck_{\cc_1}, \dots, \ck_{\cc_b}$ by construction of the $\ck_{\cc_\gb}$, 
and their diameter is $< \gd_i \le \gd_{i+1}$ by hypothesis.
The contact cubes $\ck_{\ga_1}, \dots, \ck_{\ga_m}, \ck_{\cc_1}, \dots, \ck_{\cc_b}$ therefore cover 
$\ck_1, \dots, \ck_a, \cc_1, \dots, \cc_b$ and have diameter $< \gd_{i+1}$,
that is, they have all the properties required in Claim~$(i)$.
\proofend

\ni
{\bf 2. End of the proof}

\m
\ni
Recall that $M$ is covered by the $d+1$ sets 
$$
\cs^j \,=\, \bigcup_{k=1}^\ell \phi_k \bigl( (s_k \ed^j(d))_1 \bigr),
\quad\, j = 1, \dots, d+1.
$$
Fix $j$.
By Proposition~\ref{p:new.cubes}, the set $\cs^j$ can be covered by disjoint 
contact cubes $\ck_1, \dots, \ck_a$.
This means that there exist coordinate cubes
$K_1, \dots, K_a$ and contact charts $\phi_\ga \colon \RR^d \to M$ with
$\phi_\ga (K_\ga) = \ck_{\ga}$. Choose open cubes $U_\ga$ around $K_\ga$ such that the sets $\phi_\ga (U_\ga)$ are still disjoint.
By Corollary~\ref{c:cuboid}, the cubes $U_\ga$ are contactomorphic to $\RR^d$, and so the sets $\phi_\ga (U_\ga)$ are disjoint contact balls containing the sets $\ck_\ga$.
By Proposition~\ref{p:one}, $\ck_1 \cup \dots \cup \ck_a$ can therefore be covered by a single contact ball. Hence $M$ can be covered by $d+1$ contact balls.
\proofend

\section{Proof of Theorem~2. Part~I}  \label{s:proof2}

\ni
Let $M$ be a closed connected $3$-manifold.
As in Section~\ref{s:top} we denote by $\cat (M)$ the Lusternik--Schnirelmann category of $M$.
According to~\cite{GL,OR},
\begin{eqnarray*}
\cat (M) \,=\,
\left\{
 \begin{array}{ll} 
    2 & \text{if $\pi_1(M) = \{1\}$}, \\
    3 & \text{if $\pi_1(M)$ is free and non-trivial}, \\
    4 & \text{otherwise},
 \end{array}
\right.
\end{eqnarray*}
and $\cat (M) = \B(M)$ if and only if $M$ contains no fake cells or $\cat (M)=4$.
The fundamental group $\pi_1(M)$ is free if and only if each prime summand of $M$ is a homotopy sphere or $S^2 \times S^1$ or the non-orientable $S^2$-bundle over $S^1$, see~\cite[Chapter~5]{He-76}. 
Adding the proof of the Poincar\'e conjecture and the hypothesis that $M$ is orientable,
we obtain
\begin{eqnarray*}
\B (M) \,=\,
\left\{
 \begin{array}{ll} 
    2 & \text{if $M = S^3$}, \\
    3 & \text{if $M = {\#}_k (S^2 \times S^1)$}, \\
    4 & \text{otherwise};
 \end{array}
\right.
\end{eqnarray*}
here, ${\#}_k$ denotes the $k$-fold connected sum.
Since contact $3$-manifolds are orientable and in view of
Theorem~1 we arrive at
\begin{lemma}  \label{l:top}
Consider a closed connected contact $3$-manifold $(M,\xi)$. Then
\begin{eqnarray*}
\C (M,\xi) 
\left\{
 \begin{array}{ll} 
    \in \{ 2,3,4 \} & \text{if $M = S^3$},  \\
    \in \{ 3,4 \}   & \text{if $M = {\#}_k (S^2 \times S^1)$}, \\
    =               4 & \text{otherwise}.
 \end{array}
\right.
\end{eqnarray*}
\end{lemma}

In order to complete the proof of Theorem~2
we will argue as follows.
First, we shall show that $\C (M,\xi) = \B (M)$ for every tight contact structure.
In view of Lemma~\ref{l:top} it only remains to show that $\C (M,\xi) =3$ for overtwisted contact structures on $S^3$ and ${\#}_k (S^2 \times S^1)$.
The non-existence of overtwisted discs for overtwisted structures on $\RR_{\st}^3$
will immediately imply $\C( S^3,\xi_j) \ge 3$ for overtwisted structures on $S^3$.
The main point of the whole proof will be to show that 
overtwisted structures on $S^3$ and $S^2 \times S^1$ can be covered by $3$ contact charts.
From this we shall easily obtain the same result for overtwisted structures on 
connected sums of $S^2 \times S^1$.

\begin{proposition}  \label{p:tight=}
$\C (M,\xi) = \B(M)$
for all tight contact $3$-manifolds $(M,\xi)$.
\end{proposition}

\proof
Consider a contact chart $\phi \colon \RR^3 \to M$.
Since $\xi$ is tight, $\phi^* \xi$ is a tight contact structure on $\RR^3$.
By a theorem of Eliashberg~\cite{El-92}, there is a diffeomorphism $\psi$ of $\RR^3$ such that $\psi^* (\phi^* \xi) = \xi_{\st}$.
Then $\phi \circ \psi \colon (\RR^3, \xi_{st}) \to (M,\xi)$ is a contact chart with the same image as $\phi$.
We conclude that $\C (M,\xi) \le \B(M)$, and the proposition follows.
\proofend

\begin{proposition}  \label{p:S3}
$\C (S^3,\xi_j) \ge 3$ for all overtwisted contact structures $\xi_j$ on $S^3$.
\end{proposition}

\proof
Consider an overtwisted contact structure $\xi_j$ on $S^3$.
Arguing by contradiction we assume that $S^3 = \cb_1 \cup \cb_2$, 
where $\cb_i = \phi_i (\RR^3)$ are contactomorphic images of $\RR^3_{\st}$.
Let $D \subset (S^3,\xi_j)$ be an overtwisted disc.
We can assume that the centre $p = \phi_1(0)$ of $\cb_1$ is disjoint from $D$.
Choose open balls $B_r(0)$ and $B_R(0)$ in $\RR^3$ such that
$$
\phi_1 \bigl( B_r(0) \bigr) \cap D \,=\, \emptyset
\qquad \text{ and } \qquad
M \setminus \phi_1 \bigl( B_R(0) \bigr) \,\subset\, \cb_2 .
$$
The contact vector field $X=(x,y,2z)$ with Hamiltonian $H=2z+xy$ generates the contact isotopy 
$$
\gf_H^t (x,y,z) \,=\, \left( e^tx, e^ty, e^{2t}z \right) .
$$
Choose $T$ so large that 
$$
\gf_H^T \left( \RR^3 \setminus B_r(0) \right) \,\subset\, 
\RR^3 \setminus B_R(0) ,
$$
and let $f \colon \RR^3 \to [0,1]$ be a smooth cut-off function with
$f |_{B_R(0)} =1$ and $f |_{B_{2R}(0)} =0$.
The contact isotopy $\gf_{fH}^t$, $0 \le t \le T$, of $\RR^3$
is then compactly supported and maps $\RR^3 \setminus B_r(0)$ into $\RR^3 \setminus B_R(0)$.
As a consequence, the contactomorphism $\psi$ of $M$ defined by
\begin{eqnarray*}
\psi (m) \,=\,
\left\{
 \begin{array}{ll} 
   \bigl( \phi_1 \circ \gf_{fH}^T \circ \phi_1^{-1} \bigr) (m)  & \text{if $m \in \cb_1$},  \\
   m  & \text{if $m \notin \cb_1$},
 \end{array}
\right.
\end{eqnarray*}
maps $D$ into $\cb_2$, and so $(\phi_2^{-1} \circ \psi)(D)$ is an overtwisted disc 
in $\RR^3_{\st}$, 
which is impossible in view of Bennequin's Theorem from~\cite{Ben}.
\proofend

\begin{proposition}  \label{p:le3}
$\C (M,\xi) \le 3$
for all contact structures $\xi$ on $M = S^3$ or on $M = S^2 \times S^1$. 
\end{proposition}

\ni
The idea of the proof is given at the end of this section.
The proof is postponed to Sections~\ref{s:convex} and \ref{s:proofend}.

\b
Given two contact $3$-manifolds $(M_1,\xi_1)$ and $(M_2,\xi_2)$, the connected sum $M_1 \# M_2$ 
carries a canonical contact structure, see~\cite{Gei-2001, Gei-book,We-91} 
or Section~\ref{s:ex}.5.

\begin{proposition}  \label{p:connected}
Every contact structure $\xi$ on $\#_k (S^2 \times S^1)$ can be written as a
connected sum of contact manifolds
$$
(S^2 \times S^1,\xi_1) \#\cdots \# (S^2\times S^1,\xi_k).
$$
\end{proposition}

\proof
We distinguish the cases $\xi$ tight and $\xi$ overtwisted.

\s
\ni
{\bf $\xi$ is tight.}
Assume first that the contact structure $\xi$ on $M := \#_k (S^2 \times S^1)$
is tight. We shall show that
\begin{equation}  \label{e:ind}
(M,\xi) \,\cong\, (S^2\times S^1,\xi_0)\# \cdots \# (S^2\times S^1,\xi_0).
\end{equation}
Here $(S^2\times S^1,\xi_0)$ denotes the unique tight structure on $S^2 \times S^1$,
see~\cite{El-92,Gi-91}.

We assume first that $k=2$, that is, $M = (S^2 \times S^1) \# (S^2 \times S^1)$.
Choose a smooth separating $2$-sphere $S \subset M$.
Then each of the $2$ components of $M \setminus S$ is a copy of $S^2 \times S^1$
with a $3$-ball removed.
After perturbing $S$, we can assume that $S$ is convex 
(see Section~\ref{s:convexity} below for the definition).
Since $M$ is tight,
a neighbourhood $S \times I$ of $S$ is tight.
According to~\cite{Gi-91}, there is a unique tight structure on $S \times I$.
Giroux flexibility~\cite{Gi-91} now allows us to further perturb
$S$ so that the characteristic foliation on $S \times I$ is standard,
that is, the characteristic foliation on each $S \times \{t\}$ looks like the one 
on the boundary of the $3$-ball $B := B^3_1(0)$ in $(\RR^3,\xi_{\st})$.
This allows us to glue to both boundary spheres of $M \setminus S$
a ball $B$, so as to obtain $2$ copies of $S^2 \times S^1$.

We can now apply a theorem due to Colin, see Theorem~2.6 in \cite{Co-99}, 
which states that if the complement of a convex $2$-sphere 
in a closed contact manifold $P$ is tight, then so is $P$. 
The sphere $S$ in $S^2 \times S^1$ is convex, 
and its complement is the tight ball $B$ and one of the components of 
$M \setminus S$, which is tight. 
Therefore, the contact structures on the $2$ copies of $S^2 \times S^1$ 
are tight and therefore diffeomorphic to $\xi_0$.
By the definition of contact connected sum,
the connected sum of these $2$ copies of $(S^2 \times S^1,\xi_0)$ obtained  
by using the above balls $B$ is $(M,\xi)$.

Assume now that $M = \#_k (S^2 \times S^1)$ with $k \ge 3$. 
Arguing by induction, we assume that~\eqref{e:ind} holds true for $k-1$.
Set $M' = \#_{k-1} (S^2 \times S^1)$, so $M = M' \times (S^2 \times S^1)$.
Repeating the above argument, we see that
$$
(M,\xi) \,\cong\, (M',\xi') \# (S^2 \times S^1, \xi_0)
$$
for some tight contact structure $\xi'$ on $M'$.
By the induction hypothesis, 
$$
(M',\xi') \,\cong\, \#_{k-1}(S^2 \times S^1, \xi_0) ,
$$
and hence 
$(M,\xi) \cong \#_k(S^2 \times S^1, \xi_0)$.

\m
\ni
{\bf $\xi$ is overtwisted.} 
We show that every overtwisted contact structure $\xi_{\ot}$ on 
$M=\#_k S^2 \times S^1$ can be written as 
$$
(M,\xi_{\ot}) \cong (S^2 \times S^1,\xi_1) \# \cdots \# (S^2 \times S^1,\xi_k),
$$
where each $(S^2 \times S^1,\xi_i)$ is an overtwisted contact structure 
on $S^2\times S^1$. Note that $H^2(M;\ZZ)\cong \ZZ^k$, 
so we can interpret the first Chern class as $k$ integers,
$$
c_1(\xi_{\ot})=(n_1,\ldots, n_k) \in \ZZ^k \cong H^2(M;\ZZ).
$$
By Eliashberg's classification of overtwisted structures, 
we find overtwisted structures $\tilde \xi_i$ on $S^2\times S^1$
with
$$
c_1(\tilde \xi_i) = n_i \in \ZZ \cong H^2(S^2\times S^1;\ZZ),
\qquad i=1,\ldots,k .
$$ 
(This follows, in fact, already from Martinet's existence results~\cite{Ma}). 
Consider now
$$
(M,\tilde \xi) \,:=\, (S^2 \times S^1,\tilde \xi_1) \# \cdots \# (S^2 \times S^1,\tilde \xi_k) .
$$
Then $c_1(\tilde \xi) = c_1(\xi_{\ot})$. 
Therefore, on the $2$-skeleton of $M$ (that is, on the complement of a $3$-ball in $M$) the contact structures $\xi_{\ot}$ and $\tilde \xi$ are isotopic as plane fields, 
see for instance the description of classifying plane fields 
in Chapter~11.3 of~\cite{GoS}. 
The obstruction of extending this isotopy to $M$ lies in 
$H^3(M;\ZZ) = \ZZ$.
Using again Eliashberg's classification of overtwisted contact structures, we therefore find an overtwisted contact structure 
$\xi_e$ on $S^3$ such that the isotopy extends to an isotopy of plane fields between $\xi_{\ot}$ and the contact structure $\hat \xi$ defined by 
$$
(M,\hat\xi) = 
(S^2 \times S^1,\tilde \xi_1) \# \cdots \# (S^2 \times S^1,\tilde \xi_k)\# (S^3,\xi_e).
$$
By the uniqueness part of Eliashberg's classification of overtwisted contact structures, $\xi_{\ot}$ and $\hat\xi$ are diffeomorphic.
Set now $(S^2 \times S^1,\xi_i) := (S^2 \times S^1,\tilde \xi_i)$ 
for $i<k$ and
$(S^2 \times S^1,\xi_k) := (S^2 \times S^1,\tilde \xi_k)\# (S^3,\xi_e)$.
\proofend

\begin{proposition}  \label{p:connectedsum}
$\C \left( M_1 \# M_2, \xi_1 \# \xi_2 \right) \le
\max \bigl\{ \C \left( M_1, \xi_1 \right), \C \left( M_2, \xi_2 \right) \bigr\}$
for any two closed contact $3$-manifolds $(M_1,\xi_1)$ and $(M_2,\xi_2)$.
\end{proposition}
This result holds true in any dimension, see Theorem~\ref{t:connectedsum.d} below.

\s
\b
\ni
{\bf Idea of the proof of Proposition~\ref{p:le3}}

\begin{figure}[ht] 
 \begin{center}
  \psfrag{S1}{$S^1$}
  \psfrag{S2}{$S^2$}
  \psfrag{T}{$T$}
  \leavevmode\epsfbox{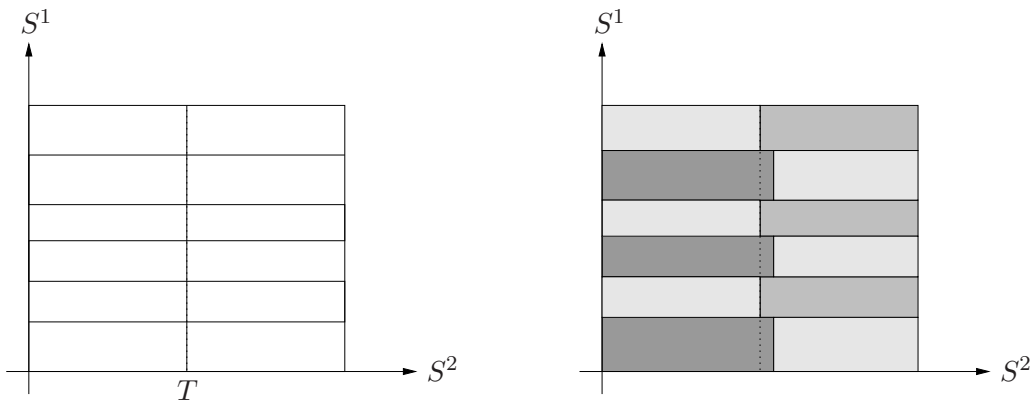}
 \end{center}
 \caption{The idea of the proof of Proposition~\ref{p:le3}.} 
 \label{figure.idea2}
\end{figure}

\m
\ni
Our aim is to show that $\C (M,\xi) \le 3$ for all contact structures on $M = S^3$ and $M = S^2 \times S^1$.
The proof goes similarly for the two cases.
We give the idea for $S^2 \times S^1$.
Consider the foliation of $S^2 \times S^1$ by the spheres $S_\tau = S^2 \times \{ \tau \}$ with $\tau \in S^1$.
The main step will be to construct a smooth $S^1$-family of embedded closed curves $\gg_\tau \subset S_\tau$ 
such that
for each $\tau \in S^1$ the curve $\gg_\tau$ divides $S_\tau$ into two closed discs with tight neighbourhoods.
Let $T = \bigcup_{\tau \in S^1} \gg_\tau$ be the torus in $S^2 \times S^1$ formed by the curves $\gg_\tau$.
This torus and a subdivision of $S^1$ yield a partition
of $S^2 \times S^1$ into pieces as shown on the left of Figure~\ref{figure.idea2}.
For a sufficiently fine partition of $S^1$, each piece has a tight neighbourhood. 
After replacing some of the pieces by slightly smaller ones and some by slightly bigger ones, 
we can paint the pieces by three colours, such that pieces of the same colour are disjoint, 
see the right of Figure~\ref{figure.idea2}. 
The claim then follows from Proposition~\ref{p:one}.
This plan will be carried out in the next two sections.

%
%

\section{Characteristic foliations and convexity}  \label{s:convex}

\ni
In this section we collect some notions and results from 
contact $3$-manifold topology that we shall use in the proof of 
Proposition~\ref{p:le3}.
The results of this section are due to Giroux~\cite{Gi-91,Gi-00}.
Throughout, $(M,\xi)$ is a contact $3$-manifold.
We assume that the contact structure $\xi$ is co-oriented
by means of a 1-form $\ga$ such that $\xi= \ker \ga$.
Let $S$ be a closed embedded oriented surface in~$M$.

\subsection{Characteristic foliations}  
\label{s:char}

An {\it oriented singular foliation}\/ $\cf$ on $S$ 
is an equivalence class of smooth vector fields on $S$,
where two vector fields $Y,Y'$ are 
equivalent if  $Y'=fY$ for a smooth positive function~$f$.
We write $\cf_Y$ for the singular foliation represented by $Y$.  
The flow lines of $Y$ are the oriented leaves of~$\cf_Y$.  

A point $p \in S$ is called a {\it singular point}\/ of $\cf_Y$ if $Y(p)=0$.
At a singular point $p$, we define the divergence $\div Y (p)$ to
be the trace of the linearisation of $Y$ at~$p$ (this agrees with
the usual definition of divergence with respect to an arbitrary 
area form on~$S$). Since $\div (fY)(p) = f \div Y(p)$, 
the sign of $\div Y(p)$ does not change when $Y$ is replaced 
with an equivalent vector field.

Pick a contact form $\ga$ co-orienting $\xi$ and an area form $\Omega$ on $S$ 
orienting~$S$.
Define the vector field $Y$ on $S$ by 
$$
i_Y \Omega \,=\,  \ga |_S ,
$$
then $Y_x$ generates $\xi_x \cap T_x S$ at the points $x$  
where $\xi$ is transverse to~$S$ and  $Y_x=0$
at the points where $\xi$ is tangent to~$S$.
Changing the contact form $\ga$ representing the co-orientation of  $\xi$
and the area form $\Omega$ 
representing the orientation of $S$ results in
multiplying $Y$ by a  positive function. 
Hence the oriented singular 
foliation~$\cf_Y$ does not depend on these choices. 
Reversing the co-orientation of $\xi$ or the  orientation of $S$ 
leads to replacing $\cf_Y$ with $\cf_{-Y}$.
We shall call $\cf_Y$ the 
{\it (oriented) characteristic foliation}\/ 
on $S$ induced by~$\xi$. 

An oriented singular foliation $\cf_Y$ on $S$ is a characteristic 
foliation of a contact structure near $S$ 
if and only if $\div Y (p) \neq 0$ at all singular points~\cite{Gi-91}.
In view of this result, we call such 
a vector field a {\it characteristic vector field},
and we call the corresponding oriented singular foliation
an oriented characteristic foliation 
even when the contact structure is not specified.
Vector fields that generate characteristic foliations 
form an open subset in the $C^\infty$-topology among all vector fields.

\subsection{Perturbations}
The following two assertions show that perturbations of
characteristic foliations can be always induced by 
perturbations of the underlying surfaces.

\begin{lemma}  \label{l:per1}
For every sufficiently $C^\infty$-small perturbation $\cf'$
of the  characteristic foliation
on $S$ induced by~$\xi$, there exists a 
$C^\infty$-small diffeomorphism $\Phi$
of $M$ such that   $\Phi$ maps $\cf'$ to the characteristic
 foliation on $\Phi (S)$.
\end{lemma}

\begin{lemma} \label{l:per2}
Let $M=S^2\times [0,1]$ resp.\ $M=S^2 \times S^1$. 
Denote by $\cf_\tau$ the  characteristic foliations
induced by $\xi$ on  $S_\tau:=S^2\times\{\tau\}$.
Assume that  $\{\cf'_\tau\}, \tau\in[0,1]$ resp.\ $\tau \in S^1$, 
is a family of oriented singular foliations
which is sufficiently $C^\infty$-close to $\{\cf_t\}$.
If $M=S^2\times [0,1]$ also assume that $\cf'_0=\cf_0$, $\cf'_1=\cf_1$. 
Then there is  a $C^\infty$-small diffeomorphism $\Phi$
of $M$ such that $\Phi$ maps $\cf'_\tau$ to the characteristic
foliation on $\Phi (S_\tau)$ induced by~$\xi$. 
\end{lemma}
\ni
Lemma~\ref{l:per1} is an obvious corollary of Lemma~\ref{l:per2}. 

\m

\proofof{Lemma~\ref{l:per2}}
Let $\ga$ be the contact form on $M = S^2 \times [0,1]$ defining $\xi$,  
let $\Omega$ be the area 
form on~$S^2$ defining orientations on all $S_\tau$,
and denote $I=[0,1]$.  
By abuse of notation, we write 
$\ga= f d\tau + \beta_\tau$, where
$f$ is a function on $S^2\times I$ and $\{\beta_\tau\}$ is a 
family of $2$-forms on $S^2$ depending on the parameter $\tau\in I$.
Consider the smooth family of vector fields $\{Y_\tau\}$, $\tau\in I$, 
on $S^2$ defined by $i_{Y_\tau}\Omega=\beta_\tau$.
For each $\tau\in I$, the lift of $Y_\tau$ to 
$S_\tau$ (by $x\mapsto (x,\tau)$)
directs the characteristic foliation $\cf_\tau$ on 
$S_\tau$ induced by~$\xi$. 
Let $\{Y'_\tau\}$, $\tau\in I$, be the family of vector fields on $S^2$
such that the lift of $Y'_\tau$ to 
$S_\tau$ generates the singular foliation  $\cf'_\tau$ for all~$\tau$,
$\{Y'_\tau\}$ is $C^\infty$-close to $\{Y_\tau\}$, and we have 
$Y'_0=Y_0$, $Y'_1=Y_1$. 
Define a 1-form $\ga'$ on $M$ by 
$\ga'=f d\tau+ i_{Y'_\tau}\Omega$.
Then $\ga'$ is $C^\infty$-close to $\ga$.
By the relative Gray's theorem, there is a
$C^\infty$-small diffeomorphism $\Phi$ of $M$
such that $\Phi^*\ga'=g\ga$, where $g$ is a positive function. 
This diffeomorphism has the required properties. 
The proof for $M=S^2 \times S^1$ is the same.
\proofend
                                      
\subsection{Singular points}
\begin{figure}[ht] 
 \begin{center}
  \psfrag{a}{(a)}
  \psfrag{b}{(b)}
  \psfrag{c}{(c)}
  \leavevmode\epsfbox{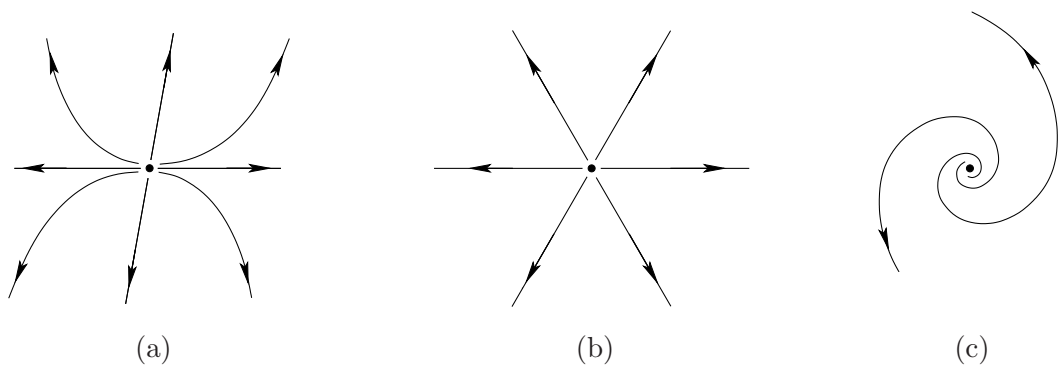}
 \end{center}
 \caption{A generic node, a bicritical node, and a focus.} 
 \label{figure.nodes}
\end{figure}

\begin{figure}[ht] 
 \begin{center}
  \psfrag{a}{(a)}
  \psfrag{b}{(b)}
  \leavevmode\epsfbox{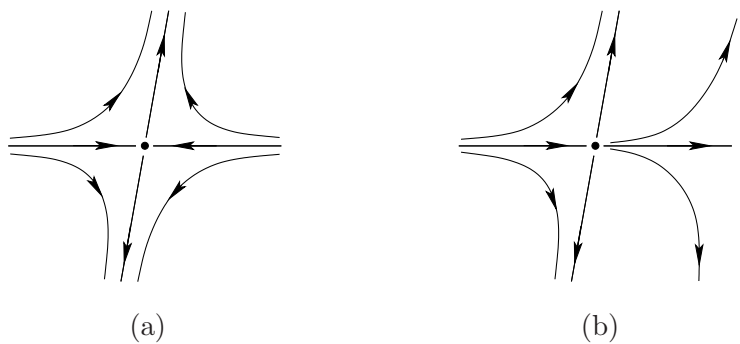}
 \end{center}
 \caption{A saddle and a saddle-node.} 
 \label{figure.saddle}
\end{figure}

We say that singular points $p \in (S,\cf_Y)$ and $p' \in (S',\cf_{Y'})$ 
have the same {\it topological type}\/
if there exists a homeomorphism between neighbourhoods $U(p) \subset S$ 
and $U(p') \subset S'$ mapping the oriented leaves of  
$\cf_Y |_{U(p)}$ to the oriented leaves of $\cf_{Y'}|_{U(p')}$.
Two examples of isolated singular points 
are a generic node 
(Figure~\ref{figure.nodes}a)
and a focus 
(Figure~\ref{figure.nodes}c). 
Both have the same topological type as the bicritical node 
(Figure~\ref{figure.nodes}b).
Figure~\ref{figure.saddle}a shows a saddle  and Figure~\ref{figure.saddle}b
a saddle-node.

\begin{lemma}  \label{l:sing}
An isolated singular point $p$ of a  characteristic 
foliation $\cf_Y$ on $S$ has the topological type 
of a node, a saddle, or a saddle-node.
\end{lemma}

\proof
Since $\div Y (p) \neq 0$, the linearisation of $Y$ at $p$ 
has at least one non-zero eigenvalue, and all its non-zero eigenvalues
have non-zero real parts. 
Then by the Shoshita\u{\i}shvili theorem~\cite{Sho0,Sho1},
there is a homeomorphism between a neighbourhood of $p$ in $S$
and a neighbourhood of the origin in $\RR^2$ that sends $p$
to the origin and maps the oriented flow lines of $Y$ to the   
oriented flow lines of the vector field 
$$
\pm x_1 \partial/\partial x_1+\,h(x_2) \,\partial/\partial x_2. 
$$
Since $p$ is an isolated singular point,
the function $h$ has an isolated zero at $x_2=0$.
If $h$ changes sign at $x_2=0$, then topologically $p$
is a node or a saddle, otherwise it is a saddle-node.   
\proofend

%
%

%
%
%
%

\ni
Let $\phi_Y^t$ denote the flow of a  characteristic vector field~$Y$.
A singular point $p$ of $Y$ is {\it positive}\/ if $\div Y (p)>0$, and  
{\it negative}\/ if $\div Y (p)<0$. 
Note that $p$ is a positive (resp.~negative) singular point
if  the orientation of $\xi_p=T_p S$ given by 
the restriction of $d\ga$ agrees 
(resp.~disagrees) with the orientation of~$S$.
A singular point $p$ is said to {\it attract}\/ (resp.~{\it repel})
an orbit $L$ of the  characteristic foliation if 
$\lim_{t \to +\infty} \phi_Y^t (x) = p$ 
(resp.~$\lim_{t \to -\infty} \phi_Y^t (x) = p$)
for $x\in L$. 
A topological node is a  {\it source}\/ (resp.~a {\it sink})
if it repels  (resp.~attracts) all orbits of $\cf_Y$ passing through 
a sufficiently small neighbourhood.
Sources are always positive and sinks are always negative. 
Similarly, there are two kinds of saddle-nodes: a saddle-source is positive
and attracts one non-constant orbit;
a saddle-sink is negative and repels one non-constant orbit.
Saddles can be positive or negative; each of them attracts two 
and repels two non-constant orbits.

\subsection{Convexity}  \label{s:convexity}

\s
\ni
A {\it retrograde connection}\/ is 
an orbit of a characteristic
foliation going from a negative singular point to a positive one. 
There are four types of retrograde connections:
saddle to saddle, saddle to saddle-node, saddle-node to saddle,
and saddle-node to saddle-node.
Note that reversing the orientation of the characteristic
foliation  leaves intact the set of retrograde connections.
%

%
%
%
\begin{figure}[ht] 
 \begin{center}
  \psfrag{-}{$-$}
  \psfrag{+}{$+$}
  \leavevmode\epsfbox{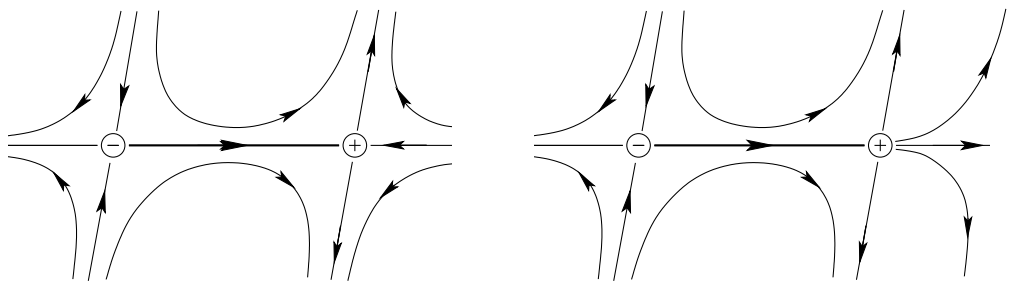}
 \end{center}
 \caption{Two retrograde connections.} 
 \label{retro}
\end{figure}
%
%

\m
\ni
A periodic orbit (or a cycle) is called {\it degenerate}\/ 
if the Poincar\'e return map has eigenvalue equal to~$1$, 
and {\it non-degenerate}\/ otherwise. 

An embedded closed orientable surface 
$S \subset (M,\xi)$ is called {\it convex}\/ if there exists a 
contact vector field $X$ on $M$ that is transverse to~$S$.
Giroux showed that convex surfaces form an open dense subset
among all embedded surfaces.
The following result is Proposition 2.5 of~\cite{Gi-00}.
\begin{proposition}  \label{p:crit.convexity}
Let $S \subset (M,\xi)$ be a $2$-sphere such 
that the characteristic foliation 
$\cf_\xi$ has only finitely many singular points. 
Then $S$ is convex if and only if $\cf_\xi$ 
has no degenerate cycles and no retrograde connections.
\end{proposition}

\subsection{A criterium for tightness}

There are simple geometric criteria that allow to find out
whether a convex surface has a tight neighbourhood. 
In order to formulate one of them, we need yet another notion.

Let $S$ be a convex surface in $(M,\xi)$. 
Given a contact vector field $X$ transverse to~$S$, 
one defines the {\it dividing set}\/ $\Delta_X\subset S$ as the set 
of points $x$ where $X(x) \in \xi_x$. 
The dividing set is a collection of disjoint
embedded circles and it does not depend, up to isotopy, on the choice
of the transverse contact vector field.
The following criterium for tightness is proved in 
\cite[Th\'eor\`eme~4.5.a)]{Gi-01}.

\begin{theorem}  \label{thm:Giroux_OTcriterion}
A convex $2$-sphere $S$ has a tight neighbourhood
if and only if
the dividing set of $S$ is connected.
\end{theorem}

There is another useful criterion, which is related to this theorem. 
Consider a characteristic foliation $\cf$ with 
isolated singular points on a $2$-sphere~$S$.
Define the {\it positive graph}\/ $\Gamma_+(\cf)$ to be the union
of the positive singular points of $\cf$ (they are the vertices 
of the graph) and the orbits of $\cf$ that connect two positive 
singular points (their closures are the edges of the graph;
the endpoints of an edge can coincide). 
Note that the number of edges of $\Gamma_+(\cf)$
is finite because each positive singular point of $\cf$ 
attracts at most two orbits.
Similarly, we define the {\it negative graph}\/  
$\Gamma_-(\cf)$ to be  the union
of the negative singular points of $\cf$  and
the orbits of $\cf$ that connect two negative singular points. 
Finally, define the {\it full graph}\/ 
$\Gamma (\cf) = \Gamma_+(\cf) \cup \Gamma_-(\cf)$.

\begin{lemma} \label{l:criterionOT}
Let $S$ be a convex $2$-sphere with 
a characteristic foliation~$\cf$. 
Then $S$ has a tight neighbourhood
if and only if
$\cf$ has no periodic orbit and its full graph $\Gamma (\cf)$ is a union of disjoint trees. 
In this case, each of the graphs $\Gamma_+(\cf)$, $\Gamma_-(\cf)$ is a tree.
\end{lemma}

\proof
Assume first that $S$ has a tight neighbourhood.
By Theorem~\ref{thm:Giroux_OTcriterion}, 
the dividing set $\Delta$ of $\cf$ is an embedded circle.
The complement $S \setminus \Delta$ consists of two open discs, $D_-$ and $D_+$.
There is a characteristic vector field~$Y$ representing~$\cf$ and an
area form $\go$ on $S \setminus \Delta$ such that, 
after possibly renaming the discs,
$\div_\go Y = -1$ on $D_-$ and $\div_\go Y = 1$ on $D_+$,
see \cite{Gi-91} or \cite[p.~230]{Gei-book}.
Arguing by contradiction, assume that~$K$ is a cycle in~$\cf$ 
or a loop in $\Gamma (\cf)$. 
Every trajectory of $Y$ intersecting $\Delta$ meets $\Delta$ transversely,
and goes from $D_+$ to $D_-$,
see again \cite{Gi-91} or \cite[p.~230]{Gei-book}.
Since $Y$ is tangent to $K$, it follows that $K$ is disjoint from
$\Delta$, say, $K \subset D_-$. 
Therefore, one of the connected components of $S \setminus K$ is contained in $D_-$. 
Then $\div Y = -1$ on this component. 
This contradicts the fact that $Y$ is tangent to $K$.

Assume now that $\cf$ has no cycles and $\Gamma (\cf)$ contains no loops.
Then all connected components of $S \setminus \Delta$ are discs,
see ``C. Fin de la d\'emonstration de la Proposition~3.1'' 
on page~658 of~\cite{Gi-91}. 
Therefore, $\Delta$ is connected, and $S$ has a tight neighbourhood
by Theorem~\ref{thm:Giroux_OTcriterion}.
The graphs $\Gamma_+(\cf)$,  $\Gamma_-(\cf)$ are trees because
each connected component of $S\setminus\Delta$
contains exactly one connected  component of $\Gamma (\cf)$.
\proofend

\section{Proof of Proposition~\ref{p:le3}}  \label{s:proofend}

\subsection{Tightening curves}

Let $S$ be an embedded $2$-sphere in a contact $3$-manifold $(M,\xi)$. 
We call a smoothly embedded circle $\gg$ in $S$
a {\it tightening curve} (with respect to $\xi$)
if each of the two closed discs in 
$S$ with boundary $\gg$  has a tight neighbourhood in~$(M,\xi)$.
The following lemma is an immediate consequence of Theorem~2.5.23 
in~\cite{Gei-book}.

\begin{lemma}  \label{l:ti1}
If two contact structures $\xi,\xi'$ induce the same 
characteristic foliation on the $2$-sphere~$S$ 
and $\gg$ is a tightening curve with respect to $\xi$, then
$\gg$ is also tightening with respect to~$\xi'$.
\end{lemma}

We can therefore say that a curve is (or is not) 
tightening with respect to a characteristic foliation $\cf$ on~$S$,
or $\cf$-tightening, 
without specifying the contact structure.
The following lemma shows that the tightening property is $C^\infty$-open.

\begin{lemma} \label{l:t-stability}
Let $\gg\subset S$ be a tightening curve with respect 
to a characteristic foliation~$\cf$. 
If  $\gg'$ is sufficiently $C^\infty$-close to $\gg$ and
$\cf'$ is sufficiently $C^\infty$-close to $\cf$, then 
$\gg'$ is tightening  with respect to~$\cf'$.
\end{lemma}

\proof 
Denote by $D_0, D_1$ the two closed discs in 
$S$ with boundary~$\gg$.  
According to Lemma~\ref{l:per1}, there exists
a $C^\infty$-small  contactomorphism $\Phi$ that  maps 
$\cf'$ to the characteristic foliation  $\cf''$ on $\Phi (S)$
induced by the contact structure. 
Denote by $D_0', D_1'$ the two closed discs in 
$S$ with boundary~$\gg'$.
If  $\gg'$ is sufficiently $C^\infty$-close to $\gg$ and
$\cf'$ is sufficiently $C^\infty$-close to $\cf$, we can assume that
$\Phi(D_j')$ is $C^\infty$-close to $D_j$, $j=0,1$.
In particular, we can achieve that the disc $\Phi(D_j')$
is contained in a tight neighbourhood of $D_j$, $j=0,1$.
Then $\Phi(\gg')$ is tightening with respect to $\cf''$
and  hence $\gg'$ is tightening with respect to~$\cf'$.
\proofend

\subsection{From families of tightening curves to contact atlases}

There are slight differences in the proof of
Proposition~\ref{p:le3} for the two cases, $M=S^3$ and $M=S^2\times S^1$.
For $M=S^2\times S^1$, consider the foliation of $M$ by the spheres 
$S_\tau=S^2\times\{\tau\}$. 

\begin{proposition}  \label{p:-S2xS1}
Assume that there exists a smooth family of oriented tightening curves 
$\gg_\tau\subset S_\tau\subset (S^2\times S^1,\xi)$, $\tau\in S^1$.
Then $C(S^2\times S^1,\xi)\le3$.
\end{proposition}

For $M=S^3$, fix two disjoint embedded closed balls, $B_0$ and $B_1$, 
in~$M$. Identify  the complement of their interiors with $S^2\times I$, 
where $I=[0,1]$, by means of a diffeomorphism. 
Consider the foliation of  $S^2\times I$ by the spheres 
$S_\tau=S^2\times\{\tau\}$.

\begin{proposition}  \label{p:-S3}
Assume that there exists a smooth family of tightening curves 
$\gg_\tau\subset S_\tau\subset (S^3,\xi)$, $\tau\in I$, 
and that each of the balls $B_0,B_1$ has a tight neighbourhood. 
Then $C(S^3,\xi)\le 3$.
\end{proposition}

\begin{figure}[ht] 
 \begin{center}
  \psfrag{S1}{$S^1$}
  \psfrag{S2}{$S^2$}
  \psfrag{D}{$D$}
  \psfrag{D'}{$D'$}
  \psfrag{D+}{$D_+$}
  \psfrag{D-'}{$D_-'$}
  \psfrag{B0}{$B_0$}
  \psfrag{B1}{$B_1$}
  \psfrag{J1}{$J_1$}
  \psfrag{J2}{$J_2$}
  \psfrag{J2k}{$J_{2k}$}
  \leavevmode\epsfbox{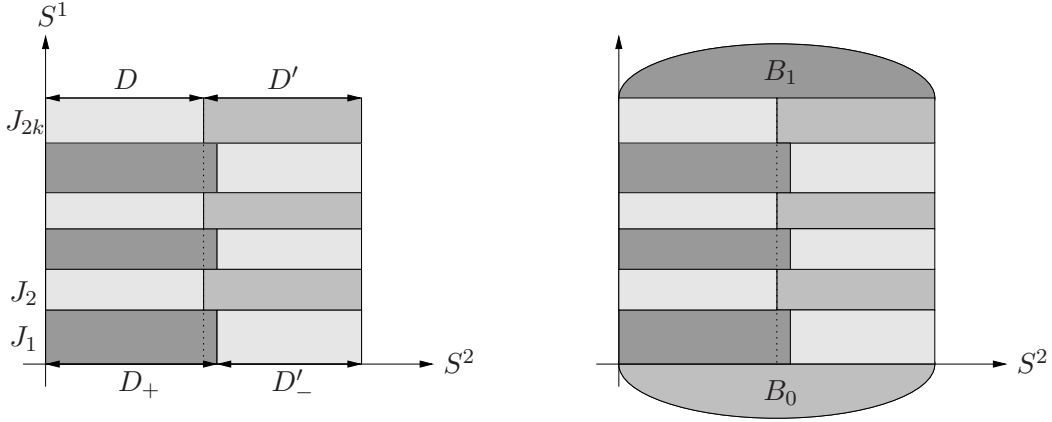}
 \end{center}
 \caption{The partition of $S^2 \times S^1$ and of $S^3$ into three sets.} 
 \label{figure.3balls}
\end{figure}

\proofof{Proposition~\ref{p:-S2xS1}} 
The union $T$ of all curves $\gg_\tau$ is diffeomorphic to the torus
due to the orientation hypothesis.
We can assume, after applying a diffeomorphism of $S^2\times S^1$
that preserves each sphere $S_\tau$ as a set, that $T=\gg\times S^1$, where
$\gg$ is a curve in $S^2$, say, the equator. 
Denote by $D,D'$ the two closed hemispheres in $S^2$ with boundary~$\gg$.
For each $\tau\in S^1$, since $\gg_\tau$ is tightening,
there exist a tight neighbourhood $U$ of $D\times\{\tau\}$ 
and a tight neighbourhood $U'$ of $D'\times\{\tau\}$.
Then there is a neighbourhood $V_\tau$ of $\tau$ in $S^1$ such that
$D\times V_\tau\subset U$ and  $D'\times V_\tau\subset U'$.
By compactness, the circle~$S^1$ can be covered by
finitely many of the neighbourhoods $V_\tau$. Subdivide $S^1$
into intervals $J_1=[a_0,a_1],J_2=[a_1,a_2],\ldots,J_{2k}=[a_{2k-1},a_0]$
such that each interval $J_i$ is covered by one of the  
neighbourhoods~$V_\tau$. 
Then each of the sets $D\times J_i$, $D'\times J_i$  
has a tight neighbourhood; denote by $W_i$ a tight neighbourhood 
of~$D\times J_i$. 
Let $\gg'$ be a parallel in the hemisphere $D'$. Denote by $D_+, D'_-$ 
the two discs with boundary  $\gg'$ such that
$D\subset D_+$, $D'_-\subset D'$.
Pick  $\gg'$ sufficiently close to  $\gg$,
then $D_+\times J_i\subset W_i$ for all~$i$. 
Define $C_1$ to be the union of the sets $D_+\times J_{2i-1}$,  
$C_2$ the union of the sets $D'\times J_{2i}$,  
$C_3$ the union of the sets $D'_-\times J_{2i-1}$ and  $D\times J_{2i}$, 
where $i\in \{1,2,\ldots,k\}$, see Figure~\ref{figure.3balls}.
Then $S^2\times S^1=C_1\cup C_2 \cup C_3$ and every connected component
of each of the sets $C_1,C_2,C_3$ has a tight neighbourhood.
For a given set $C_i$, these neighbourhoods can be chosen disjoint.
Hence, by Proposition~\ref{p:chart.3} and  
Proposition~\ref{p:one}, $C_i$ can be covered by a 
single contact ball. Thus~$C(S^2\times S^1,\xi)\le3$.
\proofend

%
%

\proofof{Proposition~\ref{p:-S3}} 
Arguing as in the proof of Proposition~\ref{p:-S2xS1},
we construct, after applying to $\xi$
a suitable diffeomorphism of  $S^2\times I$,
a subdivision of $I$ into a union of intervals
$J_1=[0,a_1],J_2=[a_1,a_2],\ldots,J_{2k}=[a_{2k},1]$
and subdivisions of $S^2$ into discs $S^2=D\cup D'$,  $S^2=D_+\cup D'_-$
such that $D\cap D'_-=\emptyset$ and each of the sets
$D\times J_i$, $D'\times J_i$, $D_+ \times J_i$, $D'_-\times J_i$, 
$i \in \{1,2,\ldots 2k\}$, 
has a tight neighbourhood.
Define $C_1$ to be the union of the sets $D_+ \times J_{2i-1}$ and of the ball $B_1$ (where $\pp B_1 = S^2 \times \{1\}$), 
$C_2$ to be the union of the sets $D' \times J_{2i}$ and of the ball $B_0$ 
(where $\partial B_0=S^2\times \{0\}$), and
$C_3$ to be the union of the sets $D'_-\times J_{2i-1}$ and  $D \times J_{2i}$
where $i\in \{1,2,\ldots,k\}$, see Figure~\ref{figure.3balls}.
Then $S^3=C_1\cup C_2 \cup C_3$ and every connected component
of each of the sets $C_1,C_2,C_3$ has a tight neighbourhood.
Hence, by Proposition~\ref{p:chart.3} and Proposition~\ref{p:one}, 
each $C_i$ can be covered by a single contact ball. 
Thus~$C(S^3,\xi)\le3$.
\proofend

In order to construct families of tightening curves involved in 
Proposition~\ref{p:-S2xS1} and Proposition~\ref{p:-S3},
it is convenient to make the characteristic foliations
on the spheres $S_\tau$ as generic as possible by a perturbation 
of the contact form. We start by recalling the necessary definitions.

\subsection{Structurally stable and quasi-generic vector fields}

A vector field $Y$ on a manifold $M$ is called {\it structurally stable} 
if for each $Y'$ sufficiently $C^1$-close to $Y$ there is a
$C^0$-small homeomorphism $g$ of $M$ that maps oriented orbits of
$Y'$ to oriented orbits of~$Y$ 
(see~\cite{AP,PdM} for a more precise formulation).
A vector field $Y$ on a $2$-dimensional sphere is structurally stable
if and only if the following conditions are satisfied~\cite{deBaggis,PdM}:

\begin{itemize}

\item[(S1)] All singular points of $Y$ are non-degenerate, 
in the sense that the eigenvalues of the linearisation 
of $Y$ at a singular point have non-zero real parts.

\item[(S2)] All periodic orbits of $Y$  are non-degenerate;
there are finitely many of them.
 
\item[(S3)] No orbit of $Y$  is a saddle-to-saddle connection.

\end{itemize}

Structurally stable vector fields on $S^2$ form an open and dense set in the
$C^\infty$-topology~\cite{Pe59,PdM}.
We say that a characteristic vector field on a $2$-dimensional 
sphere is
{\it quasi-generic}
if it satisfies the properties (S1)--(S3)  
with exactly one of the following exceptions:

\begin{itemize}

\item[(Q1)]  One of the singular points is a saddle-node; 
none of its three separatrices connects it with saddles.

\item[(Q2)]  One of the periodic orbits is degenerate; 
the second derivative of its Poincar\'e return map is non-zero.

\item[(Q3)]  There is one saddle-to-saddle connection.

\end{itemize}

This definition agrees with the classical definition for  general 
vector fields given in~\cite{So74}. 
That definition allows as a possible exception also 
one saddle with zero divergence; 
such a node cannot be a singular point of a characteristic vector field. 
The structural stability and quasi-genericity properties
do not change when the characteristic vector field is multiplied
by a positive function. 
Therefore, they extend to characteristic foliations.
We shall call a  characteristic foliation 
 (Q1)-quasi-generic or  (Q2)-quasi-generic or  (Q3)-quasi-generic
depending on which of the three exceptions is realized.

\subsection{Constructing tightening curves}

\begin{proposition}  \label{p:ti}
If the characteristic foliation $\cf$ on a $2$-sphere~$S$ is 
structurally stable or quasi-generic, then there is a tightening 
curve with respect to~$\cf$.
\end{proposition}

\ni
{\bf Remark.} Actually, the statement is also true for an arbitrary 
characteristic foliation. The proof for the general case follows the same 
approach as the one we give below but its details are more complicated.

\m

We call a smoothly embedded circle $\gg\subset S$ 
{\it extensive}\/ with respect to a characteristic foliation $\cf$
(or  $\cf$-extensive) if the following conditions are satisfied:
 
\begin{itemize}
\item[(E1)] The curve $\gg$ intersects every periodic orbit 
of~$\cf$ and the intersection is transverse at at least one point.
\item[(E2)] The curve $\gg$ intersects at least one non-constant 
orbit in each loop of the graph~$\Gamma(\cf)$
and the intersection is transverse at at least one point.
\item[(E3)] The curve $\gg$ intersects every 
saddle-to-saddle connection of~$\cf$ 
and the intersection is transverse at at least one point.
\end{itemize}

\ni
A loop in this definition is a subset homeomorphic to~$S^1$. 

\begin{proposition}  \label{p:ext}
Let $\cf$ be a structurally stable or quasi-generic
characteristic foliation on a $2$-sphere. 
If $\gg$ is extensive with respect to $\cf$, 
then $\gg$ is tightening with respect to~$\cf$.
\end{proposition}

\ni
Extensive curves obviously exist for characteristic foliations
that are  structurally stable or quasi-generic. 
Hence Proposition~\ref{p:ti}
follows from Proposition~\ref{p:ext}.

\m \ni

\proofof{Proposition~\ref{p:ext}}
The circle $\gg$ divides $S$ into two closed discs.
Take one of them and call it~$D$.
We are to show that $D$ has a tight neighbourhood in~$(M,\xi)$. 
Choose a smoothly embedded circle $\gg^+$ in $S$ disjoint from $D$ 
and so close to $\gg$ that $\gg^+$ is also extensive.
Denote by $D_+ \supset D$ the closed disc in $S$ bounded by $\gg^+$.
Our plan is to construct a characteristic foliation 
$\widehat \cf$ on $S$ that
satisfies the assumptions of Lemma~\ref{l:criterionOT} and
coincides with $\cf$ on~$D_+$. 

We claim that if such an $\widehat \cf$ exists, then $D$
has a tight neighbourhood.
Indeed, consider a contact structure $\hat \xi$ on a
neighbourhood of $S$ that induces the  
characteristic foliation~$\widehat \cf$.
By Lemma~\ref{l:criterionOT}, there is a  neighbourhood $U$
of $S$ such that the restriction of  $\hat \xi$ to $U$ is tight.
It follows from Theorem~2.5.23 in~\cite{Gei-book} that
there exist open neighbourhoods $V_1, V_2$ of $\Int D_+$ in 
$M$ and a diffeomorphism $V_1 \to V_2$ that acts as the identity on 
$\Int D_+$ and maps $\hat \xi$ to~$\xi$.
Since $\bigl( U \cap V_1, \hat \xi \bigr)$ is tight, we conclude that 
$\bigl(U \cap V_2, \xi\bigr)$ is a tight neighbourhood of~$D$.

\begin{lemma}  \label{l:Qi}
Let $\cf$ be a {\rm(Q2)}-quasi-generic or a  {\rm(Q3)}-quasi-generic
characteristic foliation on a $2$-sphere $S$, and let
$\gg^+$ be an $\cf$-extensive curve that
bounds a disc $D_+$. Then there is a  $C^\infty$-small perturbation~$\cf'$
of $\cf$ coinciding with $\cf$ on $D_+$ such that $\cf'$ is
structurally stable and $\gg^+$ is $\cf'$-extensive. 
\end{lemma}
\proof
Assume that  $\cf$ is  {\rm(Q2)}-quasi-generic. Let $K$ be the degenerate 
periodic orbit of~$\cf$. Let $W$ be a neighbourhood of $K$ in $S$
such that $\cf$ is transverse to~$\partial W$.
Since $\gg^+$ is extensive, the set  $K\setminus D_+$ is nonempty.
Pick a point $x\in K\setminus D_+$. 
Let $U$ be a neighbourhood of $x$ which is contained in $W$ and disjoint from~$D_+$. 

\begin{figure}[ht] 
 \begin{center}
  \psfrag{K}{$K$}
  \psfrag{K1}{$K_1$}
  \psfrag{K2}{$K_2$}
  \leavevmode\epsfbox{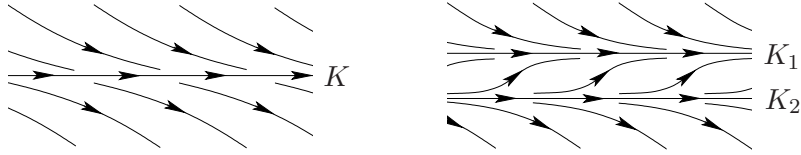}
 \end{center}
 \caption{Breaking the degenerate cycle $K$ into two non-degenerate cycles $K_1,K_2$.} 
 \label{figure.K}
\end{figure}

By a $C^\infty$-small perturbation of  $\cf$ with support in $U$, we construct a 
characteristic foliation  $\cf'$ such that  $\cf'|_W$ 
has exactly two periodic orbits $K_1$ and $K_2$, 
both non-degenerate, and no singular points. 
Figure~\ref{figure.K} illustrates schematically the effect of this perturbation 
on the characteristic foliation.
It follows from the fact that $\cf'$ is transverse to $\partial W$
that each orbit intersecting  $\partial W$ has 
one of the cycles $K_1,K_2$
as its $\alpha$- or $\omega$-limit set. Then 
$\cf'$ has no saddle-to-saddle connections passing through $W$ and hence 
no saddle-to-saddle connections at all. 
Therefore, $\cf'$ is  structurally stable.

%
%

We show now that $\gg^+$ is $\cf'$-extensive.
Condition~(E1) holds for $K_1,K_2$ because they are $C^\infty$-close
to~$K$, and for all other cycles because they were not modified. 
Since an orbit passing through $\partial W$ has 
$K_1$ or $K_2$ as its $\alpha$- or $\omega$-limit set, it cannot
be an edge of the graph~$\Gamma(\cf')$.
Thus $\Gamma(\cf)=\Gamma(\cf')$ and (E2) holds.
Finally, (E3) is satisfied because $\cf'$  is structurally stable.

Assume that  $\cf$ is  {\rm(Q3)}-quasi-generic.
Note that by (S1), the singular points of $\cf$ are nodes and saddles.
Suppose first that the saddle-to-saddle connection $L$ is heteroclinic,
that is, it goes from a saddle $p$ to a different saddle~$q$.
Then $\cf$ has no polycycles, and by the Poincar\'e--Bendixson Theorem,
every $\ga$- and $\omega$-limit set of an orbit is either a node, a saddle, or a cycle.  
Pick a point $x\in L\setminus D_+$ (it exists since $\gg^+$ is extensive).
There exists a neighbourhood $U$ of $x$  that is disjoint from~$D_+$,
such that for each point $y$ in $U\setminus L$ the orbit of $\cf$ passing 
through  $y$ has the same $\alpha$-limit set as one of the two incoming
separatrices of $p$ and the same  $\omega$-limit set as one of the two 
outgoing separatrices of~$q$. Let $\cf'$ be a generic perturbation of $\cf$
with support in~$U$. The  perturbation being generic,
we can assume that $\cf'$ has  no orbit that goes from  $p$ to $q$
and coincides with $L$ outside~$U$. 
Denote by $L_p$ (resp.~$L_q$)  
the orbit of $\cf'$ that coincides with $L$ near $p$ (resp.~$q$).
By our choice of~$U$, 
each orbit of $\cf'$ that passes through $U$ and differs from $L_p$ 
has a node or a limit cycle as its $\alpha$-limit set.
Hence such an orbit cannot be a cycle or a saddle-to-saddle-connection.
By our choice of  $\cf'$,
the trajectory $L_p$ differs from $L_q$ 
 and hence it has  a node or a limit cycle as its $\omega$-limit set.
Thus $\cf'$ has no saddle-to-saddle connections or cycles passing through~$U$,
and hence is structurally stable. 

The curve $\gg^+$ is $\cf'$-extensive because we have
$\Gamma(\cf')\subset\Gamma(\cf)\cup \overline{L_p} \cup \overline{L_q}$, and
$\gg^+$ transversely intersects at least once each of $L_p$ and $L_q$
provided that $\cf'$ is sufficiently $C^\infty$-close to~$L$.

\begin{figure}[ht] 
 \begin{center}
  \psfrag{B}{$B$}
  \psfrag{Bs}{$B_*$}
  \psfrag{p}{$p$}
  \psfrag{K}{$K$}
  \psfrag{x}{$x$}
  \psfrag{U}{$U$}
  \psfrag{W}{$W$}
  \psfrag{L}{$L$}
  \psfrag{L}{$L$}
  \psfrag{Lo}{$L_\mathrm{o}$}
  \psfrag{Li}{$L_\mathrm{i}$}
  \psfrag{Lo'}{$L_\mathrm{o}'$}
  \psfrag{Li'}{$L_\mathrm{i}'$}
  \leavevmode\epsfbox{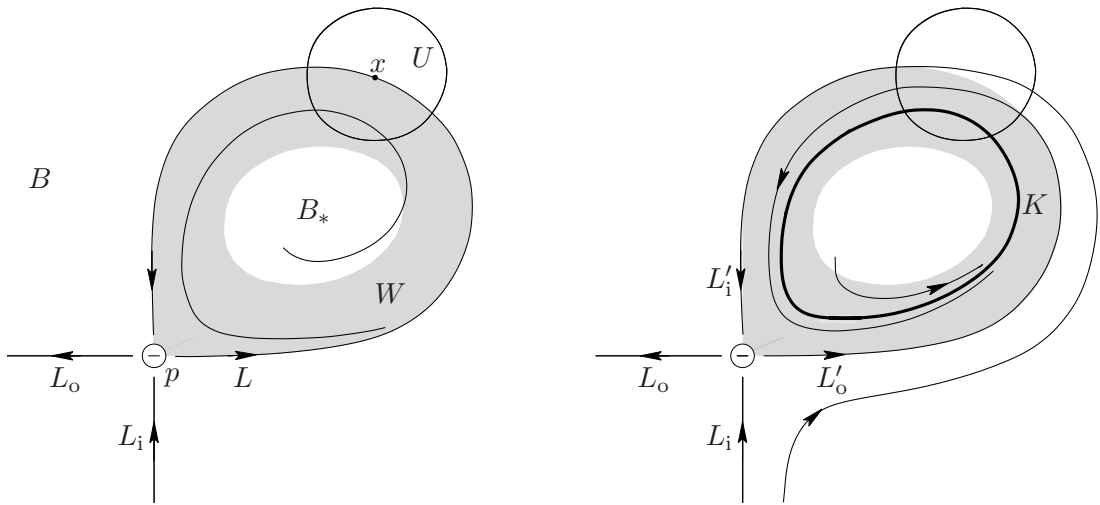}
 \end{center}
 \caption{Perturbing a homoclinic cycle into a cycle.}
 \label{figure.limitcycle}
\end{figure}

Assume now that the saddle-to-saddle connection $L$ is homoclinic,
that is, it connects a saddle~$p$ with itself.
Every $\ga$- and $\omega$-limit set of an orbit is a node, a saddle, a cycle, or the polycycle $L \cup \{p\}$.  
Denote by $L_\mathrm{o} $ (resp.~$L_\mathrm{i}$) the 
outgoing  (resp.~incoming) separatrix of $p$ different from~$L$.
Denote by $B$ the connected component of $S\setminus (L\cup\{p\})$
that contains  $L_\mathrm{o}$ and~$L_\mathrm{i}$,
and by $B_*$ the one that does not.
Assume for definiteness that the saddle $p$ is negative.
According to Theorem~44 in \S29 of \cite{ALGM} the polycycle $L \cup \{p\}$ is attracting from one side, that is, there exists an open set $W \subset B_*$, 
which is the intersection of a neighbourhood of $L \cup \{p\}$ with $B_*$, 
such that the $\omega$-limit set of each $y \in W$ is $L \cup \{p\}$.
Pick a point $x\in L\setminus D_+$. 
There is a neighbourhood $U$ of $x$  disjoint from~$D_+$, 
$L_\mathrm{o}$, $L_\mathrm{i}$ such that $U \cap B_* \subset W$ and  
such that each orbit of $\cf$ that passes
through $U\cap B$ has the same $\alpha$-limit set  as $L_\mathrm{i}$
and the same $\omega$-limit set as~$L_\mathrm{o}$,
see the left of Figure~\ref{figure.limitcycle}.

%
%

As illustrated on the right of Figure~\ref{figure.limitcycle},
there exists an arbitrarily $C^\infty$-small perturbation
$\cf'$ of  $\cf$ supported in $U$ such that
$\cf'$ has exactly one limit cycle $K$ passing through~$U$,
the cycle $K$ is non-degenerate and it is 
the $\omega$-limit set for each orbit 
of $\cf'$ entering $W$ as well as for the 
outgoing separatrix $L'_\mathrm{o}$ of $p$ different from~$L_\mathrm{o}$
(for the proof, see e.g.~Sections~2 and 3 of \S29 of~\cite{ALGM}).
Denote by $L'_\mathrm{i}$ the incoming  separatrix of $p$ for $\cf'$
different from~$L_\mathrm{i}$.

The curve $\gg^+$ is $\cf'$-extensive because we have
$\Gamma(\cf')\subset\Gamma(\cf)\cup \overline{L'_\mathrm{o}} 
\cup \overline{L'_\mathrm{i}}$, and
$\gg^+$ transversely intersects at least once each of 
$L_\mathrm{o}'$, $L_\mathrm{i}'$, and $K$ (near the point where $\gg^+$ 
transversely intersects~$L$)
provided that $\cf'$ is sufficiently $C^\infty$-close to~$L$.

Suppose that  $L'$ is a saddle-to-saddle
orbit of  $\cf'$. Then it has to pass  through~$U$.
Let $y\subset \partial U$ be the point of $\overline{U}$
that $L'$ passes last. If $y\in B$, then the $\omega$-limit of $L'$
is the same as for~$L_\mathrm{o}$ and hence not a saddle,
a contradiction. 
If $y\in B_*$, then the $\omega$-limit of $L'$ is $K$,
a contradiction.
If  $y\in L$, then $L'=L'_\mathrm{i}$, and the $\go$-limit set of $L'$ is the one of $L_i$, and hence not a saddle.
This completes the proof of Lemma~\ref{l:Qi}.
\proofend

We continue the proof of Proposition~\ref{p:ext}.
By Lemma~\ref{l:Qi}, we can assume that $\cf$
is structurally stable or {\rm(Q1)}-quasi-generic.
We now modify $\cf$ outside~$D_+$,
not restricting to small perturbations anymore,
with the goal to eliminate the periodic orbits of~$\cf$.

\begin{figure}[ht] 
 \begin{center}
  \leavevmode\epsfbox{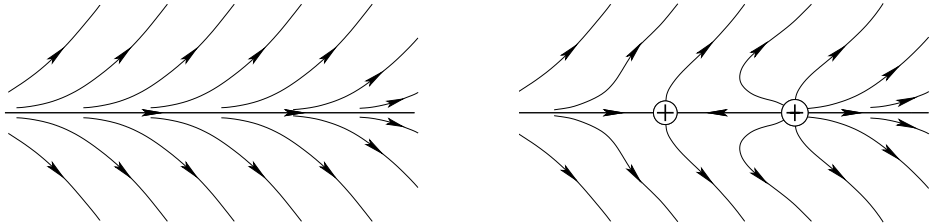}
 \end{center}
 \caption{Breaking a limit cycle.}
 \label{figure.break}
\end{figure}

\begin{figure}[ht] 
 \begin{center}
  \psfrag{q}{$q$}
  \psfrag{r}{$r$}
  \leavevmode\epsfbox{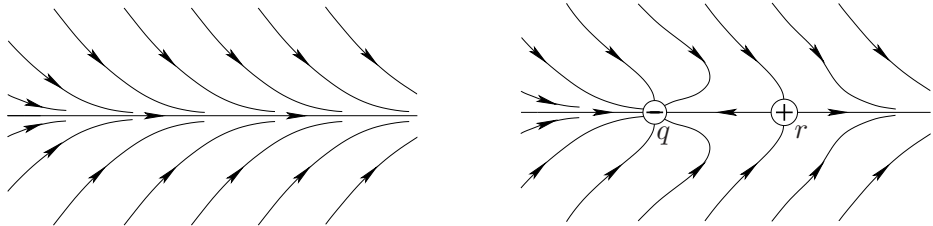}
 \end{center}
 \caption{How to eliminate a loop in $\Gamma_-(\cf)$.}
 \label{figure.qr}
\end{figure}

Let $K$ be a (necessarily non-degenerate) limit cycle of~$\cf$.
There exists a foliation of  a neighbourhood $W$ of $K$ in $S$
into circles $K_s$, $s\in \left]-1,1\right[$, such that $K_0=K$
and $\cf$ is  transverse to all $K_s$ with $s\ne0$.
Pick a point $x\in K\setminus D_+$ and a neighbourhood $U\subset W$ of $x$
disjoint from $D_+$. 
We replace the characteristic foliation $\cf$ with
a characteristic foliation $\cf'$ that coincides with $\cf$ 
outside~$U$, is transverse to all $K_s$ with $s\ne0$, 
and has exactly two singular points with nonzero divergence,
a saddle and a node. 
More precisely, for a repelling limit cycle we insert
a positive saddle and a source as shown in Figure~\ref{figure.break},
and for an attracting limit cycle we insert a sink and a negative saddle, cf.~Figure~\ref{figure.qr}. 
This operation does not create new cycles. 
Since there are only finitely many separatrices of $\cf$
that enter or leave $U$, for a generic choice of a foliation $\cf'$
with the properties described above, $\cf'$ 
has no saddle-to-saddle connections. 
Thus  $\cf'$ is   structurally stable or {\rm(Q1)}-quasi-generic.
The curve $\gg^+$ is $\cf'$-extensive because we have
$\Gamma(\cf')=\Gamma(\cf)\cup K$ and  $\gg^+$ intersects
$K$ transversely.

%
%

Applying this procedure to all cycles of $\cf$ in succession,
we construct a characteristic foliation such that it coincides with  $\cf$
on $D_+$, it is  structurally stable or {\rm(Q1)}-quasi-generic,
it has no cycles, and $\gg^+$ is $\cf'$-extensive with respect to it.
We can thus assume that $\cf$ has these properties.  
Note that $S$ is convex in view of Proposition~\ref{p:crit.convexity}.

At the next step of our construction, we eliminate loops in the 
graph~$\Gamma(\cf)$. If there are no loops, that is,
$\rk H_1(\Gamma(\cf))=0$, then by Lemma~\ref{l:criterionOT} the $2$-sphere~$S$ 
has a tight neighbourhood and Proposition~\ref{p:ext}
is proved. 
Assume that $\rk H_1(\Gamma(\cf))>0$.  
Assume for definiteness that $\rk H_1(\Gamma_-(\cf))>0$.
Let $P$ be a loop in~$\Gamma_-(\cf)$.
Pick a non-singular point $x$ in $P\setminus D_+$. 
Denote by $L$ the orbit of $\cf$ passing through~$x$.
This orbit arrives in a singular point $p$,
which is either a sink or a saddle-sink,
and in the latter case $L$  arrives at the sink side of~$p$
(that is, $L$ is not a parabolic separatrix).
There exists a neighbourhood $U$ of $x$  with the
following properties: 
\,(1)\,  $U$ is disjoint from~$D_+$;
\,(2)\,  $U\cap \Gamma(\cf)= U\cap L$;
\,(3)\,  each trajectory of $\cf$ passing through $U$ has
connected intersection with~$U$;
\,(4)\, each trajectory of $\cf$ passing through $U$ arrives in~$p$,
and if $p$ is a saddle-sink, then it arrives at the sink side of~$p$.

We construct a characteristic foliation $\cf'$ that coincides
with $\cf$ on $D_+$ and has two singular points in $U$,
a sink $q$ and a positive saddle~$r$, such that the outgoing
separatrices of $r$ arrive in $p$ and $q$ (see Figure~\ref{figure.qr}).
Since only finitely many separatrices of $\cf$ enter $U$,
we may choose $\cf'$ in such a generic way that each of the two 
incoming separatrices of $r$
comes from a source or from the source side of a saddle-source.

%
%

We claim that $\cf'$ is  structurally stable or {\rm(Q1)}-quasi-generic,
that $\rk H_1(\Gamma(\cf'))=\rk H_1(\Gamma(\cf))-1$,
 and that  $\gg^+$ is $\cf'$-extensive.
The first claim follows from the fact that every trajectory of $\cf'$
leaving $U$ arrives in~$p$, and every trajectory of $\cf'$ that
enters $U$ and stays inside $U$ arrives either in $q$ or in~$r$.
Let $L_\mathrm{i},L'_\mathrm{i}$ denote
 the incoming separatrices of $r$,
and let $L'$ denote the orbit of $\cf'$ that 
arrives in $q$ and coincides with the ``negative half'' of $L$
outside~$U$. Then we have
\[
\Gamma_+(\cf')=\Gamma_+(\cf)\sqcup L_\mathrm{i}\sqcup \{r\}\sqcup
L'_\mathrm{i}, \qquad
\Gamma_-(\cf')=\bigl(\Gamma_-(\cf)\setminus L\bigr) \sqcup L'\sqcup \{q\}.
\]
It follows immediately that 
$\rk H_1(\Gamma_-(\cf'))=\rk H_1(\Gamma_-(\cf))-1$.
The newly added piece of $\Gamma_+(\cf')$ connects two
vertices in   $\Gamma_+(\cf)$ that belong to different 
connected components of $S\setminus P$, and hence
to different connected components of $\Gamma_+(\cf)$.
Thus $\rk H_1(\Gamma_+(\cf'))=\rk H_1(\Gamma_+(\cf))$
and $\rk H_1(\Gamma(\cf'))=\rk H_1(\Gamma(\cf))-1$.
The curve  $\gg^+$ is $\cf'$-extensive because every
loop in $\Gamma(\cf')$ is a loop in~$\Gamma(\cf)$.

Iterating this loop elimination procedure, 
we produce a  characteristic foliation without loops.
This completes the proof of Proposition~\ref{p:ext}.
\proofend

\subsection{Constructing families of tightening curves}
For $M=S^2\times S^1$, consider the foliation of $M$ by the spheres 
$S_\tau=S^2\times\{\tau\}$.
Denote by $\cf_\tau$ the characteristic foliation induced on 
$S_\tau$ by~$\xi$. 
It follows from Theorem 2 in~\cite{So74} that one can $C^\infty$-approximate
the family of characteristic foliations $\{\cf_\tau\}$
by a family  of characteristic foliations $\{\cf_\tau'\}$
where each  $\cf'_\tau$ is either structurally stable or 
quasi-generic and the values of $\tau$ for which $\cf'_\tau$ 
is structurally stable form an open and dense subset $Z\subset S^1$. 
By Lemma~\ref{l:per2}, 
there is a diffeomorphism $\Phi$ of $S^2\times S^1$ 
such that  $\cf'_\tau$ is induced by $\Phi^*\xi$ for each~$\tau$.
After replacing $\xi$ with $\Phi^*\xi$, 
we can assume that  $\cf_\tau$ is  structurally stable when $\tau\in Z$
and $\cf_\tau$ is quasi-generic when $\tau\in S^1\setminus Z$.

By Proposition~\ref{p:ti} there exists for each $\tau \in S^1$
an $\cf_\tau$-tightening curve~$\gamma'_\tau \subset S_\tau$.  
By Lemma~\ref{l:t-stability} there exists for each $\tau\in S^1$
a neighbourhood $V_\tau$ of $\tau$ in $S^1$ such that
the translate of $\gamma'_{\tau}$ into $S_{\tau'}$ 
is an $\cf_{\tau'}$-tightening curve for every $\tau'\in V_\tau$. 
By the compactness of the circle  $S^1$, it can be covered by
finitely many of the neighbourhoods~$V_\tau$. 
Subdivide $S^1$
into intervals $J_1=[a_0,a_1],J_2=[a_1,a_2],\ldots,J_{2k}=[a_{2k-1},a_0]$,
such that each  $J_i$ is covered by some neighbourhood  $V_{\tau_i}$
and all the endpoints $a_i$ belong to~$Z$. 
Given $\tau'\in J_i$, denote by  $\gamma^i_{\tau'}$ the translate 
of  $\gamma'_{\tau_i}$ into~$S_{\tau'}$. 
The curve $\gamma^i_{\tau'}$ is $\cf_{\tau'}$-tightening.

\begin{lemma}  \label{l:interpolation}
Let $S$ be a $2$-sphere with a structurally stable
characteristic foliation $\cf$
and let $\gamma, \gamma'$ be  oriented
 $\cf$-tightening curves.
Then    $\gamma$ and $\gamma'$ can be connected by a smooth
path in the space of  oriented $\cf$-tightening curves.
\end{lemma}

\ni
The proof of this lemma is postponed until the end of this section.

Denote by $\cll$ the union over all $\tau\in S^1$ of
the spaces of oriented smoothly embedded circles in~$S_\tau$,
with the $C^\infty$-topology.
Denote by $\pi$ the natural projection $\cll\to S^1$ that
sends $\gamma\in S_\tau$ to~$\tau$. 
We construct a piecewise-smooth map $\psi\colon S^1 \to \cll$ as follows. 
Divide $S^1$ into $4k$ intervals, $I_1,\ldots,I_{4k}$. 
For each $i\in\{1,2,\ldots,2k\}$, 
let $\gs_i \colon I_{2i-1} \to J_i$ be an orientation preserving diffeomorphism,  
and map $\tau \in I_{2i-1}$ to the curve $\gg^i_{\gs_i(\tau)}$.
We equip these curves with an orientation.
The interval $I_{2i}$ is mapped to a family 
of oriented  $\cf_{a_{i}}$-tightening curves in $S_{a_{i}}$
that connects $\gamma^{i}_{a_{i}}$  to $\gamma^{i+1}_{a_{i}}$;
such a family exists by Lemma~\ref{l:interpolation}.
Then $\pi\circ\psi$ maps $I_{2i-1}$ to $J_i$ and $I_{2i}$ to~$a_i$.

By Lemma~\ref{l:t-stability}, for each  $s\in S^1$ there exists
a neighbourhood $U_s$ of $\psi(s)$ in $\cll$ such that
each curve $\gamma\in U_s$ is  an  $\cf_\tau$-tightening curve 
in~$S_\tau$, where $\tau=\pi(\psi(s))$. 
The union $U$ of the sets $U_s$ over all $s\in S^1$ is
a neighbourhood of the set $\psi(S^1)$ in~$\cll$.
There exists  a smooth map $\gf\colon S^1\to \cll$  such that
$\gf(S^1)\subset U$ and $\pi\circ\gf$ is a diffeomorphism from 
$S^1$ to~$S^1$.
Given $\tau \in S^1$, define
$\gamma_\tau=\gf(s)\subset S_\tau$, where $\pi(\gf(s))=\tau$.
The curve  $\gamma_\tau$ is $\cf_\tau$-tightening for each $\tau \in S^1$.
Hence  $C(S^2\times S^1,\xi)\le3$ by Proposition~\ref{p:-S2xS1}.

The proof for the case $M=S^3$ goes as follows.
Pick two points $p_0, p_1 \in S^3$, and choose 
local coordinates $(x, y, z)$ on disjoint neighbourhoods $U_0\supset p_0$,
 $U_1\supset p_1$ such that
$$
\xi= \ker (dz+xdy-ydx),\,\,\, p_0,p_1=(0,0,0).
$$  
For $\eps>0$ small enough, the balls 
$B_0,B_1=\{x^2+y^2+z^2\le \eps\}$ are contained in
the tight neighbourhoods~$U_0, U_1$. 
We identify the complement of the interiors of these balls
with~$S^2\times I$.
The characteristic foliation induced by $\xi$
on each of the spheres $\partial B_0, \partial B_1$ 
is a singular foliation with two singular points, the poles $x=y=0$.
The non-constant leaves of this singular foliation connect the poles.
It is structurally stable. 
This allows us to $C^\infty$-approximate the family of characteristic foliations $\{\cf_\tau\}$, $\tau \in I$, by a family of structurally stable or quasi-generic characteristic foliations $\{ \cf_\tau' \}$ in such a way that $\cf_0 = \cf_0'$ and $\cf_1=\cf_1'$.
Arguing as in the case $M=S^2\times S^1$,
we construct, using  Lemma~\ref{l:interpolation},
a smooth family formed by tightening curves
$\gamma_\tau\subset S_\tau=S^2\times\{\tau\}$, $\tau\in I$.
Then $C(S^3,\xi)\le3$ by Proposition~\ref{p:-S3}.

\subsection{Proof of Lemma~\ref{l:interpolation}}
In view of Lemma~\ref{l:t-stability}, we can assume,
after a perturbation, that the tightening curves $\gamma,\gamma'$
are transverse to all cycles in $\cf$ and all orbits that are parts 
of the graph~$\Gamma(\cf)$.
We then claim that the curves   $\gamma,\gamma'$  intersect
each  cycle in $\cf$ and each loop in~$\Gamma(\cf)$.
Since $\cf$ is structurally stable,  there are no saddle-to-saddle
connections and hence this claim implies $\gamma,\gamma'$
to be $\cf$-extensive. 

Suppose the claim fails, say, for~$\gamma$.
Denote by $K$ the cycle in $\cf$ or the loop in $\Gamma(\cf)$
not intersected by~$\gamma$. 
Denote by $D$ the disc in $S$ that is
bounded by $\gamma$ and contains~$K$.
Let $U$ be a tight neighbourhood of $D$. 
Pick a $2$-sphere $S'\subset U$ that contains~$D$.
There is a  $C^\infty$-small perturbation $S''\subset U$
of $S'$ such that the characteristic foliation $\cf''$ induced
on $S''$ by the contact structure is structurally stable.
It follows from the structural stability of $\cf$ that
$\cf''$ also has a cycle or a loop in $\Gamma(\cf'')$,
which is $C^0$-close to~$K$, provided that $S''$ is sufficiently close to $S'$.
Then, by Lemma~\ref{l:criterionOT}, $S''$ has no tight
neighbourhood, a contradiction.

Denote by $X$ (resp.~$X'$) the set of the points where 
$\gamma$ (resp.~$\gamma'$) intersects a cycle of~$\cf$
or an orbit belonging to the graph~$\Gamma(\cf)$.
We can assume that $\gamma$ intersects $\gamma'$
and the intersection is transverse. 
Indeed, otherwise we pick a small closed piece $J$ 
of $\gamma$ disjoint from~$X$ and deform $\gamma$
by an isotopy of $S$ supported outside $\gamma\setminus J$
to a curve $\gamma^*$ that intersects $\gamma'$ transversely at some point.
This deformation goes through curves that are extensive, 
and hence tightening.
Therefore, we can replace $\gamma$ by a generic 
$C^\infty$-small perturbation of~$\gamma^*$, which
intersects transversely $\gamma'$, all cycles in $\cf$, 
and all orbits that are parts of the graph~$\Gamma(\cf)$.

There exists a disjoint collection  $\{P_1,\ldots,P_m\}$,
where each $P_i$ is a piece of an orbit of $\cf$
diffeomorphic to a closed interval, and each point in $X \cup X'$
is an interior point of one of~$P_i$.   
Since the curves $\gamma,\gamma'$ have nonempty transverse intersection,
there exist two points $q_1,q_2\in S^2$ disjoint from  $\gamma,\gamma'$
such that the oriented curves $\gamma$ and $\gamma'$
represent the same nontrivial element in~$\pi_1(S^2\setminus\{q_1,q_2\})$. 
We can choose the points $q_1, q_2$ also disjoint from the intervals~$P_i$.  

\begin{figure}[ht] 
 \begin{center}
  \psfrag{g}{$\gamma$}
  \psfrag{g'}{$\gamma'$}
  \psfrag{P1}{$P_1$}
  \psfrag{Pm}{$P_m$}
  \psfrag{S}{$S^1 \times \{0\}$}
  \leavevmode\epsfbox{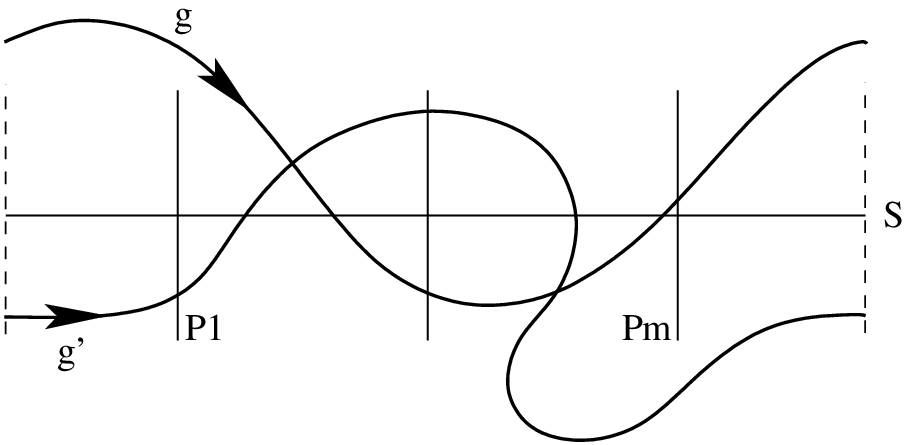}
 \end{center}
 \caption{}
 \label{figure.bi}
\end{figure}

By means of a diffeomorphism,
we can identify $S\setminus \{q_1,q_2\} $ with
the cylinder $S^1\times \left]-2; 2\right[$
in such a way that each $P_i$ is identified with the set 
$\{b_i\}\times  \left[-1; 1\right]$ for some~$b_i\in S^1$, cf.~Figure~\ref{figure.bi}.

%
%

Then each of the curves $\gamma,\gamma'$ is isotopic to the
equator~$S^1\times\{0\}$.
We perform such an isotopy in two steps.
At the first step, the curve is squeezed into
$S^1\times \left]-1; 1\right[$ by compressing along 
the second coordinate. At the second step,
the compressed curve is isotoped to the equator
inside $S^1\times \left]-1; 1\right[$
in such a generic way that at each moment
it intersects each meridian $P_i=\{b_i\}\times  \left[-1; 1\right]$
transversely at at least one point.
The same property is automatically satisfied for the curves involved
in the first step of deformation.
Therefore, this deformation goes through curves that  are extensive, 
and hence tightening.
Concatenating the paths connecting $\gamma$ and  $\gamma'$ with the
equator, we construct a path through tightening curves that
connects $\gamma $ with~$\gamma'$.
The orientations of $\gamma$ and $\gamma'$ extend to the
interpolating family of curves since $\gamma$ and $\gamma'$ 
represent the same element in~$\pi_1(S^2\setminus\{q_1,q_2\})$. 
This completes the proof of  Lemma~\ref{l:interpolation}.
Proposition~\ref{p:le3} and Theorem~2 are therefore also proved. 
\proofend


\section{A few results in higher dimensions}  \label{s:ex}

\ni
The only closed connected $1$-manifold is the circle $S^1$, 
and $\C (S^1,\xi) =2$ for the unique (and trivial) contact structure.
The contact covering numbers $\C (M, \xi)$ of closed $3$-manifolds are given by Theorem~2.
Not too much is known about the existence of contact structures on manifolds of dimension~$\ge 5$, 
see however \cite{Bou,El-91,Gei-97.1,Gei-97.2,Gei-2001,Gi-02,Us}.
In this section we look at contact manifolds of arbitrary dimension 
and prove a few results on the contact covering numbers for some special classes 
of such manifolds.

\b
\ni
{\bf 1. Spaces of co-oriented contact elements} (cf.~\cite[Example~3.45]{MS-98} and \cite[Section~1.5]{EKP-06}{\bf .} 
Consider a smooth connected manifold $N$. 
Let $S$ be a hypersurface in the cotangent bundle $T^*N$ which is fibrewise star-shaped with respect to the zero section.
This means that the fibrewise radial vector field $p\;\! \pp_p$ on $T^*N$ 
is transverse to~$S$.
The $1$-form $p\;\!dq = \sum_i p_i\:\!dq_i$ restricts to a contact form on $S$;
indeed, $p\;\!dq = i_{p \:\!\pp_p} (dp \wedge dq)$ 
for the symplectic form $dp \wedge dq$ on $T^*N$.
Denote by $\xi_S$ the corresponding contact structure on~$S$.
Given another fibrewise star-shaped hypersurface~$S'$, the contact manifolds
$(S,\xi_S)$ and $(S',\xi_{S'})$ are contactomorphic via projection along the 
vector field~$p\;\!\pp_p$.
The equivalence class of such contact manifolds is called the spherisation of $N$ and is denoted by $(S^*N,\xi)$.

\begin{proposition}  \label{p:PTN}
$\B (S^*N) \le \C (S^*N, \xi) \le 2 \min \{ \B (N), \dim N \}$.
\end{proposition} 
\begin{remarks}
{\rm
(i)
For $N=S^2$ we have $S^*N = \RP^3$, and so both inequalities are equalities.

\s
(ii)
{\it If $N$ is orientable and has vanishing Euler characteristic, then $\cl (S^*N) = \cl (N)+1$, 
and so $\cl (N)+2 \le \B (S^*N)$.}

\m
\ni
\proof
Set $n = \dim N$.
Assume first that $n = 2$.
Then $N$ is the $2$-torus and $S^*N$ is the $3$-torus, and 
$\cl (S^*N) = 3 = \cl(N) +1$.
Assume now that $n \ge 3$.
Set $k = \cl (N)$.
Assume first that $k=1$.
Then Poincar\'e duality implies that $N$ is a homology sphere.
Therefore $\cl (S^*N) = 2$.
Assume now that $k \ge 2$.
Let $p \colon S^*N \to N$ be the projection.
Since the Euler characteristic of~$N$ vanishes, $p$ admits a section.
Let $a \in H^{n-1}(S^*N)$ be its Poincar\'e dual.
The Leray--Hirsch Theorem~\cite[Theorem~4D.1]{Ha} 
asserts that each element of $H^*(S^*N)$
can be uniquely written as 
$$
p^* b + p^* b' \cup a, \quad \text{ where $b,b' \in H^*(N)$}.
$$
Choosing $b_1, \dots, b_k$ with $\prod_{i=1}^k b_i \neq 0$ in $H^*(N)$,
we see that 
$$
\prod_{i=1}^k p^* b_i \cup a \,=\, p^* \left(\prod_{i=1}^k b_i \right) \cup a 
\,\neq\, 0.
$$
Therefore, $\cl (S^*N) \ge k+1 = \cl (N) +1$.
To prove the converse, we argue by contradiction and assume that $\cl (S^*N) =: \ell \ge k+2$.
Choose $d_1, \dots, d_\ell \in \tilde H^*(S^*N)$ with $\prod_{i=1}^\ell d_i \neq 0$.
By the Leray--Hirsch Theorem, we can write
$d_i = p^* b_i + p^* b_i' \cup a$
where $\deg b_i \ge 1$ and $\deg b_i' \ge 0$.
Since $n \ge 3$ we have $\deg (a \cup a \cup a) = 3n-3 > 2n-1$, whence $a \cup a \cup a = 0$.
We can therefore compute
\begin{eqnarray*}
\prod_{i=1}^\ell d_i 
&=& 
\prod_{i=1}^\ell (p^* b_i+p^* b_i' \cup a) \\
&=& 
p^* \left( \prod_{i=1}^\ell b_i \right) +
p^* \left( \sum_{i=1}^\ell (\pm) b_i' \cup \prod_{j \neq i}^\ell b_j \right) \cup a +
p^* \left( \sum_{i_1 \neq i_2}^\ell (\pm) b_{i_1}' \cup b_{i_2}' \cup 
\prod_{j \neq i_1,i_2}^\ell b_j \right) \cup a \cup a .
\end{eqnarray*}
Since $\ell - 1 \ge k+1 > \cl (N)$, the terms $\prod_{i=1}^\ell b_i$ and $\prod_{j \neq i}^\ell b_j$ vanish.
Moreover, $k \ge 2$ implies that $\ell \ge k+2 \ge 4$. Therefore,
$$
\deg \left( \left( \prod_{j \neq i_1,i_2}^\ell b_j \right) \cup a \cup a \right) \,\ge\, 
2+n+n \,>\, 2n+1 \,=\, \dim S^*N ,
$$ 
and so these terms vanish also.
Therefore, $\prod_{i=1}^k d_i = 0$, a contradiction.
\proofend
}
\end{remarks}

\ni
{\it Proof of Proposition~\ref{p:PTN}.}
In view of Theorem~1 we only need to prove $\C (S^*N, \xi) \le 2 \B (N)$.
Let $\beta \colon \RR^n \to N$ be a smooth chart. The embedding   
\[
T^*\RR^n \to T^*N, \quad (q,p) \mapsto 
\left( \beta(q), \left(\left[ d \beta (q)\right]^T\right)^{-1}p \ \right)
\]
preserves the $1$-form $p\,dq$, and hence restricts to a contact embedding
$S^*\RR^n \to  S^*N$.
It remains to show that $S^*\RR^n \cong \RR^n \times S^{n-1}$ can be covered by $2$~contact charts.

Given an $m$-dimensional manifold~$L$, its $1$-jet space is the $(2m+1)$-dimensional manifold $J^1 L = \RR \times T^*L$. Its canonical contact structure $\xi_{\jet}$ is the kernel of the $1$-form $du - \gl$, where $u \in \RR$ and $\gl$ is the $1$-form $P\;\!dQ$ on $T^*L$.
A diffeomorphism $\gf \colon L \to L'$ between manifolds yields a contactomorphism 
\[
J^1L \to J^1L', \quad (u,Q,P) \mapsto 
\left( u, \gf(Q), \left(\left[ d \gf (Q)\right]^T\right)^{-1}P \ \right).
\] 
Therefore, $\left(J^1\RR^m,\xi_{\jet}\right)$ is contactomorphic to $\left(J^1 (S^m \setminus \{p\}),\xi_{\jet}\right)$ for any point $p \in S^m$.
Note that the linear diffeomorphism $J^1 \RR^m \to \RR^{2m+1}$,
$$
(u,Q,P) \,\mapsto\, \bigl( z(u,Q,P), x(u,Q,P), y(u,Q,P) \bigr) \,:=\, (u,-P,Q)
$$
is a contactomorphism between $\left( J^1 \RR^m,\xi_{\jet}\right)$ and $\RR^{2m+1}_{\st}$.
It follows that $\C (J^1 S^m, \xi_{\jet}) = 2$.
Proposition~\ref{p:PTN} now follows from 

\begin{lemma}
$(J^1 S^{n-1},\xi_{\jet})$ is contactomorphic to $(S^*\RR^n,\xi)$.
\end{lemma}

\proof
Let $S^{n-1}$ be the unit sphere in $\RR^n$, and denote by $\| \!\cdot \!\|$ and $\langle \cdot , \cdot \rangle$ the Euclidean norm and scalar product in $\RR^n$.
We identify $S^*\RR^n$ with $\RR^n \times S^{n-1}$ 
and $T^*S^{n-1}$ with
$\left\{ (Q,P) \in \RR^{n-1} \times \RR^{n-1} \mid \| Q \|=1, \langle Q,P \rangle =0 \right\}$.
The map $\psi \colon \RR^n \times S^{n-1} \to \RR \times T^*S^{n-1}$ defined by
$$
\psi (q,p) \,=\, \bigl( u(q,p), Q(q,p), P(q,p) \bigr)
\,:=\, \bigl( \langle q,p \rangle, p, q-\langle q,p \rangle p \bigr)
$$
is a diffeomorphism.
Moreover, the $1$-form $p\;\!dq = d (\frac 12 \|p\|^2)$ vanishes on $S^{n-1}$, and hence
\begin{eqnarray*}
\psi^* (du - P\;\!dQ) 
&=& d \langle q,p \rangle - \bigl( q-\langle q,p \rangle p \bigr) dp \\
&=& p \;\! dq + q \;\! dp - q \;\! dp \\
&=& p \;\! dq ,
\end{eqnarray*}
that is $\psi$ is a contactomorphism.
\proofend

\ni
{\bf 2. Products with surfaces.}
A contact structure $\xi$ on $M$ is said to be {\it co-orientable}\/
if $\xi$ is the kernel of a globally defined $1$-form. 
We consider a closed manifold $M$ of dimension $2n-1$ with a 
co-orientable contact structure, 
and a closed oriented surface $\Sigma$ of genus~$\ge 1$.
According to~\cite{Bou}, the product $M \times \Sigma$ carries a contact structure. 
If $\cl (M) = 2n-1$, 
then $\cl \left( M \times \Sigma \right) = 2n+1$, 
and so~\eqref{e:EE} and Theorem~1 imply that
\[
\C (M \times \Sigma,\xi) \,=\, 2n+2 
            \quad \text{\it for every contact structure $\xi$ on $M \times \Sigma$.}
\] 
In particular, $\C (T^{2n+1},\xi) = 2n+2$ for all contact structures on tori.

\b
\ni
{\bf 3. Quotients of homotopy spheres.}
Assume that $M$ is non-trivially covered by a homotopy sphere.
Then $\cat (M) = \dim M + 1$ by a result of Krasnoselski, see~\cite{GL},
and so $\C(M,\xi) = \dim M+1$ for all contact structures on $M$.

\b
\ni
{\bf 4. Higher dimensional spheres.}
The standard contact structure $\xi_0$ on the unit sphere $S^{2n+1} \subset \RR^{2n+2}$ is given by the contact form
\[
\ga_0 \,=\, \sum_{j=1}^{n+1} (x_j dy_j - y_j dx_j)
\]
where $(x_1,y_1,\dots, x_{n+1},y_{n+1})$ are Cartesian coordinates on $\RR^{2n+2}$.
Since for every point $p \in S^{2n+1}$ the manifold 
$(S^{2n+1} \setminus \{p\}, \xi_0)$ is contactomorphic to $(\RR^{2n+1}, \xi_{\st})$, 
see~\cite[Section~2.1]{Gei-handbook}, 
we have $\C(S^{2n+1},\xi_0)=2$.

\begin{proposition}  \label{p:2b}
Assume that $\xi$ is a contact structure on $S^{2n+1}$ 
such that $\C (S^{2n+1}, \xi) =2$.
Then the complement $(S^{2n+1} \setminus \{p\}, \xi)$ of any point 
$p \in S^{2n+1}$
is contactomorphic to a subset of $\RR^{2n+1}_{\st}$.
\end{proposition} 

\proof
We write again $d = 2n+1$.
After applying a diffeomorphism to~$(S^d, \xi)$, 
we can assume that $\xi = \xi_0$ on an open neighbourhood $\cn$ of~$p$.
Let $\pi \colon S^d \setminus \{p\} \to \RR^d$ be the stereographic projection.
There exists a diffeomorphism $\rho$ of $\RR^d$ such that $\rho (0) =0$ and
$\rho_* \pi_* \xi_0 = \xi_{\st}$, see~\cite[Section~2.1]{Gei-book}.
Set $\overline \pi = \rho \circ \pi$.

\begin{figure}[ht] 
 \begin{center}
  \psfrag{pi}{$\overline \pi$}
  \psfrag{S}{$S^d$}
  \psfrag{R}{$\RR^d$}
  \psfrag{p}{$p$}
  \psfrag{0}{$0$}
  \psfrag{aa}{$\cd$}
  \psfrag{a'}{$\cn$}
  \psfrag{d}{$\widehat D$}
  \psfrag{d'}{$\widehat N$}
  \psfrag{bb1}{$\cb_1$}
  \psfrag{b1}{$\widehat B_1$}
  \psfrag{b2}{$B_2$}
  \leavevmode\epsfbox{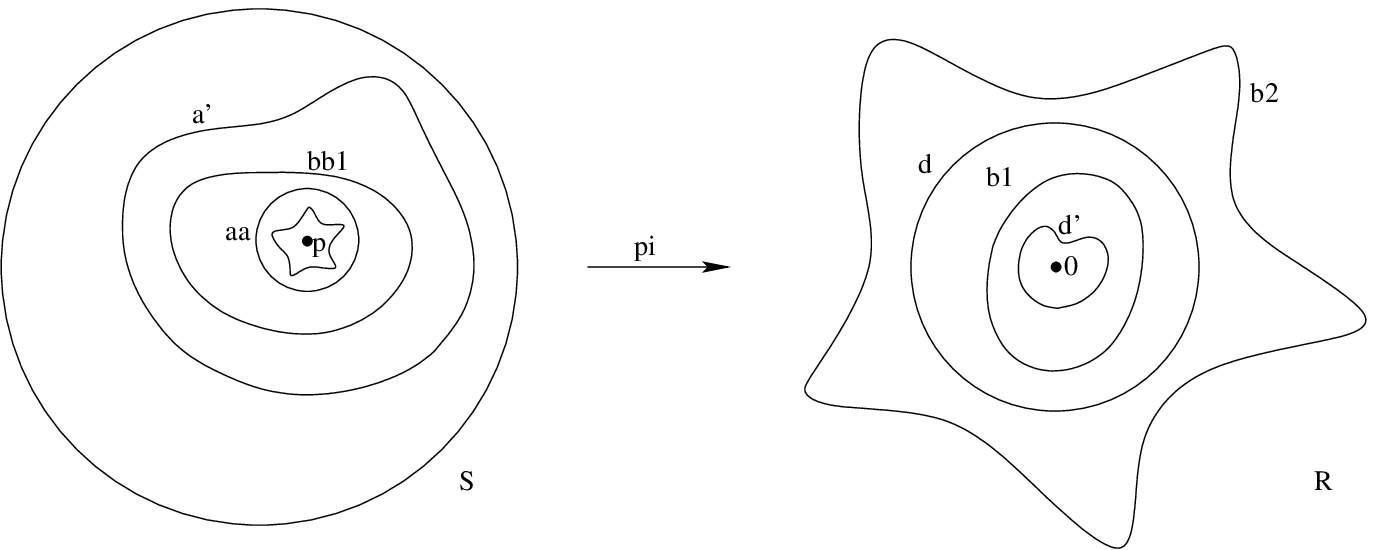}
 \end{center}
 \caption{}
 \label{figure.sets}
\end{figure}

\begin{lemma}
There is a covering  of $(S^d, \xi)$ by
contact balls $\cb_1,\cb_2$,
and an open ball $\widehat D$ in~$\RR^d$ centred at the origin
such that, denoting  $\cd := S^d \setminus \overline \pi^{-1}(\widehat D)$, we have 
$$
S^d \setminus \cb_2 \,\subset\, \cd \,\subset\, \cb_1 \,\subset\, \cn ,
$$
see the left of Figure~\ref{figure.sets}.
\end{lemma}

\proof
Since $\C(S^d,\xi)=2$, there exist contact charts 
$\phi_1,\phi_2\colon (\RR^d,\xi_{\st}) \to (S^d, \xi)$
such that $S^d$ is covered by the contact balls 
$\cb'_1 = \phi_1(\RR^d)$ and $\cb'_2 = \phi_2(\RR^d)$.
Without loss of generality
we can assume that $p \notin \cb'_2$ and that $p = \phi_1 (0)$.
Choose an open ball $E$ in $\RR^d$ centred at the origin
and so large that the contact balls
$\cb''_1 := \phi_1 (E)$ and $\cb'_2$ 
still cover~$S^d$.
Recall from Section~\ref{s:charts} that the contact Hamiltonian 
$H(\xx,\yy,z) = 2z+\xx \yy$ on $\RR^d$ generates the 
contact dilations 
\begin{equation}  \label{e:flow}
(\xx, \yy, z) \,\mapsto\, \left( e^t \xx, e^t \yy, e^{2t} z \right)
\end{equation} 
of $(\RR^d, \xi_{\st})$.
Let $f_1 \colon \cb'_1 \to [0,1]$ be a smooth 
compactly supported function with 
$f_1 |_{\cb''_1} =1$.
Let $\Phi^t$ be the contact flow on $(S^d,\xi)$ generated by the contact 
Hamiltonian $f_1 (H \circ \phi_1^{-1})$.
Choose $T_1 > 0$ so large that 
$\Phi^{-T_1} \left( \cb_1'' \right) \subset \cn$. 
 Then $S^d$ is covered by the contact balls
$\cb_1 := \Phi^{-T_1} \left( \cb_1'' \right)$, 
$\cb''_2= \Phi^{-T_1} \left( \cb'_2 \right)$, and we have
$S^d \setminus \cb''_2 \subset \cb_1 \subset \cn$.
Next, choose an open ball $\widehat D$ in $\RR^d$  centred at the origin
and so large that the set 
$\cd := S^d \setminus \overline \pi^{-1} (\widehat D)$ is contained in~$\cb_1$.
Finally,
let $f_2 \colon \cb_1 \to [0,1]$ be a smooth compactly supported 
function with 
$f_2 |_{\cb_1} =1$.
Let $\Psi^t$ be the contact flow on $(S^d,\xi)$ generated by the contact 
Hamiltonian $f_2 (H \circ \phi_1^{-1})$.
Choose $T_2 > 0$ so large that 
$ \Psi^{-T_2} \left( \cb''_2 \right) \subset S^d \setminus \cd$. 
Then $S^d$ is covered by the contact balls 
 $\cb_1$, $\cb_2:=\Psi^{-T_2} \left( \cb''_2 \right)$, 
and we have $S^d \setminus \cb_2 \subset \cd \subset \cb_1 \subset \cn$.
\proofend

%
%

The images $\widehat N := \RR^d \setminus \overline \pi (\cn)$, 
$\widehat B_1 := \RR^d \setminus \overline \pi (\cb_1)$, 
$\widehat D = \RR^d \setminus \overline \pi(\cd)$, $B_2 := \overline \pi (\cb_2)$  
in $\RR^d$ look as in Figure~\ref{figure.sets}.
Since $\xi = \xi_0$ on $\cn$, we have $\overline \pi_* \xi= \xi_{\st}$ 
on $\RR^d \setminus \widehat N$.
Let $f \colon \RR^d \to [0,1]$ be a smooth function with 
$f |_{\widehat N} =0$ and $f |_{\RR^d \setminus \widehat D} =1$.
Then $X_{fH}$ is a contact vector field on $\bigl( \RR^d, \overline \pi_* \xi\bigr)$ that makes
$\bigl( \widehat D, \overline \pi_* \xi\bigr)$ contact star-shaped.  
Proceeding exactly as in the proof of Proposition~\ref{p:contact.starshaped},
we find a contactomorphism 
$(\widehat D , \overline \pi_* \xi) \xrightarrow{\psi} 
( \RR^d, \overline \pi_* \xi\bigr)$.
A contact embedding 
$\bigl( S^d \setminus \{p\}, \xi \bigr) \ha \RR^d_{\st}$ 
is now obtained by the 
composition of contactomorphisms and an inclusion
$$
\bigl( S^d \setminus \{p\}, \xi \bigr)  \,\xrightarrow{\overline \pi}\,  
\bigl( \RR^d , \overline \pi_* \xi \bigr)  \,\xrightarrow{\psi^{-1}}\,
\bigl( \widehat D, \overline \pi_* \xi \bigr)  \,\subset\,
\bigl( B_2 , \overline \pi_* \xi \bigr)  \,\xrightarrow{(\overline \pi)^{-1}}\, 
\bigl( \cb_2 , \xi \bigr)  \,\xrightarrow{\phi_2^{-1}}\, 
\RR^d_{\st}.
$$
This completes the proof of Proposition~\ref{p:2b}.
\proofend

The concept of an overtwisted disc in a contact $3$-manifold has been 
generalized to higher dimensions in~\cite{N-06},
leading to a definition of overtwistedness in all dimensions.
It has been shown in~\cite{NvK-07} that every sphere $S^{2n+1}$ carries an overtwisted contact structure.
It has been proved in~\cite{N-06} that $(\RR^{2n+1},\xi_0)$ is not overtwisted.
Together with Proposition~\ref{p:2b} we obtain

\begin{proposition}  \label{p:overtwisted.spheres}
Let $\xi$ be an overtwisted contact structure on $S^{2n+1}$.
Then $\C (S^{2n+1}, \xi) \ge 3$.
\end{proposition} 

\ni
{\bf Remark.}
We do not know whether this result holds also true for 
$\widetilde \C (S^{2n+1}, \xi)$ if $n \ge 2$.

\m
It follows that for overtwisted contact structures on spheres,
$\B (S^{2n+1}) < \C (S^{2n+1}, \xi)$.
This shows that the contact invariant $\C(M,\xi)$ can be bigger than 
the smooth invariant $\B(M)$ in every dimension.
Problem~9.5 posed by Lutz in~\cite{Lu-88} is, however, still open in dimension $\ge 5$:

\b
\ni
{\bf Question.}
Is it true that $\C (S^{2n+1},\xi_0)=2$ if and only if $\xi = \xi_0$?

\b
\ni
Indeed, for $n \ge 2$ there are contact structures on $S^{2n+1}$ which are neither 
standard nor overtwisted, see~\cite{El-91,Mo-76,Us}.

\b
\ni
{\bf 5. Connected sums.}
The aim of this paragraph is to prove 

\begin{theorem}  \label{t:connectedsum.d}
$\C \left( M_1 \# M_2, \xi_1 \# \xi_2 \right) \le
\max \bigl\{ \C \left( M_1, \xi_1 \right), \C \left( M_2, \xi_2 \right) \bigr\}$
for any two closed contact manifolds $(M_1,\xi_1)$ and $(M_2,\xi_2)$ of the same dimension.
\end{theorem} 

\ni
{\bf Construction of the contact connected sum.}
We start with giving a precise construction of the connected sum of two contact manifolds,
which follows closely the construction of the connected sum of two smooth manifolds.
For a different description see~\cite{Gei-2001, Gei-book,We-91}.
We shall be using the rotationally symmetric contact form
$$
\alpha_{\rot} (\xx,\yy,z) \,:=\, dz+ \xx \;\! d \yy - \yy \:\! d \xx \,=\, d z + \sum_{i=1}^n r_i^2 \;\!d \phi_i
$$
on $\RR^{2n+1}$.
Here, $\xx, \yy \in \RR^n$ and $z \in \RR$,  
and $(r_1,\phi_1, \dots, r_n,\phi_n) = (\rr, \pphi)$ are multi-polar coordinates on $\RR^{2n}$.
For the linear diffeomorphism
\begin{equation}  \label{e:cont}
\psi (\xx,\yy,z) \,=\, (\xx,\yy, 2z+\xx \yy)
\end{equation}
we have $\psi^* \ga_{\rot} = 2 \ga_{\st}$, and so the contact structure $\xi_{\rot} = \ker \ga_{\rot}$ on $\RR^{2n+1}$
is contactomorphic to $\xi_{\st}$.
Note that the vector field $V(\xx,\yy, z) = (\xx,\yy,2z)$ is still a contact vector field for $\xi_{\rot}$, since
$\cll_V \ga_{\rot} = 2 \ga_{\rot}$.
Consider the unit sphere 
$$
S^{2n} \,=\, \left\{ (\xx,\yy,z) \in \RR^{2n+1} \mid \| \xx \|^2 + \| \yy \|^2 +z^2 =1 \right\}.
$$
For $t \in \RR$ set $S_t = \phi_V^t (S^{2n})$.
\begin{lemma}  \label{l:shells}
There is a contactomorphism $\Psi$ of $\left( \RR^{2n+1} \setminus \{0\}, \xi_{\rot} \right)$ 
that maps $S_t$ to $S_{-t}$ for all $t \in \RR$.
\end{lemma}

\proof
We identify $\RR^{2n+1} \setminus \{0\}$ with $S^{2n} \times \RR$ via the diffeomorphism
$$
S^{2n} \times \RR \to \RR^{2n+1} \setminus \{0\}, \qquad \mu (\ss,t) = \phi_V^t(\ss) .
$$
Since $\cll_V \ga_{\rot} = 2 \ga_{\rot}$, we have
$$
\mu^* \ga_{\rot} (\ss,t) \,=\, e^{2t} \bigl( \gb(\ss) + f(\ss) \;\! dt \bigr) ,
$$
where $\gb = \ga_{\rot} |_{S^{2n}}$ and $f(\ss) = \ga_{\rot}(\ss) (V(\ss)) = 2z$.
Define the diffeomorphism $\rho$ of $S^{2n}$ by
$$
\rho (\rr, \phi_1, \dots, \phi_n, z) \,=\, (\rr, -\phi_1 -2z, \dots, -\phi_n - 2z, z) .
$$
Then 
\begin{equation}  \label{e:mmh}
\rho^* \gb = - \gb \quad \text{ and } \quad \rho^* f = f.
\end{equation}
The diffeomorphism $\psi$ of $S^{2n} \times \RR$ defined by 
$$
\psi (\ss, t) \,=\, \left( \rho (\ss), -t \right)
$$
maps $S^{2n} \times \{t\}$ to $S^{2n} \times \{-t\}$ for all $t \in \RR$.
Moreover, by \eqref{e:mmh},
$$
\psi^* \bigl( \gb(\ss) + f(\ss) \;\!dt \bigr) \,=\, 
\rho^* \gb (\ss) - \rho^* f (\ss) \;\!dt \,=\, 
-\gb(\ss) - f(\ss) \;\!dt ,
$$
and hence $\psi$ is a contactomorphism of $\left( S^{2n} \times \RR, \mu_* \xi_{\rot} \right)$.
The diffeomorphism $\Psi := \mu^{-1} \circ \psi \circ \mu$ of $\RR^{2n+1} \setminus \{0\}$ is as required. 
\proofend

Recall that a contact structure $\xi$ is said to be 
{\it co-orientable}\/
if $\xi$ is the kernel of a globally defined $1$-form. 
A {\it co-orientation} of $\xi$ is the choice of such a $1$-form
up to multiplication by a positive function.
Consider now two contact manifolds $(M_i,\xi_i)$, $i=1,2$, 
of dimension~$2n+1$.
If $\xi_i$ is co-orientable, we assume that a co-orientation is fixed.
For $i=1,2$ choose contact charts $\phi_i \colon \left( \RR^{2n+1},\xi_{\rot} \right) \to (M_i,\xi_i)$ that preserve the existing co-orientations.
For $t \in \RR$ we set $B_t = \phi_V^t \left( B^{2n+1} \right)$, 
where $B^{2n+1}$ is the open unit ball in $\RR^{2n+1}$ centred at the origin.
The boundary of $B_t$ is~$S_t$.
Let~$\Psi$ be the diffeomorphism from Lemma~\ref{l:shells}.
The diffeomorphism 
$$
\phi_2 \circ \Psi \circ \phi_1^{-1} \colon 
\phi_1 \left(\RR^{2n+1} \right) \to \phi_2 \left( \RR^{2n+1} \right)
$$
restricts to a diffeomorphism 
$\Phi \colon \phi_1 \left(B_1 \setminus \overline{B_{-1}}\right) \to \phi_2 (B_1 \setminus \overline{B_{-1}})$.
Let $M_1 \# M_2$ be the smooth manifold obtained from 
$\bigl( M_1 \setminus \overline{\phi_1 (B_{-1})} \bigr) \cup \bigl( M_2 \setminus \overline{\phi_2 (B_{-1})} \bigr)$ 
by identifying $\phi_1 (B_1 \setminus \overline{B_{-1}})$ and $\phi_2 (B_1 \setminus \overline{B_{-1}})$ 
via~$\Phi$.
Since $\Psi$ is a contactomorphism, we can
define a contact structure $\xi_1 \# \xi_2$ on $M_1 \# M_2$ by
\begin{eqnarray*}
\xi_1 \# \xi_2 \,=\,
\left\{
 \begin{array}{ll} 
   \xi_1  & \text{on $M_1 \setminus \phi_1(B_{-1})$},  \\
   \xi_2  & \text{on $M_2 \setminus \phi_2(B_{-1})$} .
 \end{array}
\right.
\end{eqnarray*}
The contact structure $\xi_1 \# \xi_2$ on $M_1 \# M_2$ is unique up to contact isotopy.
Indeed, it follows from the Contact Disc Theorem \cite[Theorem~2.6.7]{Gei-book} 
that all contact charts in a non-co-orientable contact manifold are isotopic 
and that all co-orientation preserving contact charts in a co-oriented 
contact manifold are isotopic.

For later use we define the ``infinite neck'' 
\begin{equation}  \label{def:neck}
\cn \,:=\, \phi_1 \left( \RR^{2n+1} \setminus B_{-1} \right) \cup \phi_2 \left( \RR^{2n+1} \setminus B_{-1} \right)
\end{equation}
in $M_1 \# M_2$. 
The map $\phi \colon \left( \RR^{2n+1} \setminus \{0\}, \xi_{\rot} \right) \to \left( \cn, \xi_1 \# \xi_2 \right)$ defined by
\begin{eqnarray}  \label{def:phi}
\phi (p) \,=\,
\left\{
 \begin{array}{rl} 
   (\phi_1 \circ \Psi) (p)  & \text{if \,$p \in B_1 \setminus \{0\}$},  \\
   \phi_2(p)                & \text{if \,$p \in \RR^{2n+1} \setminus B_{-1}$},
 \end{array}
\right.
\end{eqnarray}
is a contactomorphism.

\b
\ni
{\bf Proof of Theorem~\ref{t:connectedsum.d}.}
Let now $k_1 = \C(M_1,\xi_1)$ and $k_2 = \C(M_2,\xi_2)$.
Choose contact charts 
$\phi^j_1 \colon \RR^{2n+1} \to \cb^j_1 \subset M_1$, $j = 1, \dots, k_1$, that cover $M_1$ and                        
contact charts 
$\phi^j_2 \colon \RR^{2n+1} \to \cb^j_2 \subset M_2$, $j = 1, \dots, k_2$, that cover $M_2$.
Since $k_1$ is minimal, there exists a point $p \in \RR^{2n+1}$ such that
$\phi_1^1(p)$ is disjoint from $\cb_1^2, \dots, \cb_1^{k_1}$.
After precomposing $\phi_1^1$ with an affine contactomorphism of~$\RR^{2n+1}$
that maps $0$ to $p$ (see \eqref{e:tau}), 
we can assume that $p=0$.
Recall that $B_t = \phi_V^t \left( B^{2n+1} \right)$,
where $B^{2n+1}$ is the open unit ball in $\RR^{2n+1}$ centred at the origin.
After precomposing $\phi_1^1$ with the map $\phi_V^{-t}$ for some large~$t$,
we can in fact assume that $\phi_1^1 (B_1)$ is disjoint from $\cb_1^2, \dots, \cb_1^{k_1}$.
Similarly, we can assume that $\phi_2^1 (B_1)$ is disjoint from $\cb_2^2, \dots, \cb_2^{k_2}$.
Choose $T>1$ so large that the $k_1$ contact balls $\cb^j_1(T) := \phi^j_1 \left( B_T \right)$ 
still cover $M_1$ and the $k_2$ contact balls $\cb^j_2(T) := \phi^j_2 \left( B_T \right)$ still cover $M_2$.
Form the connected sum $(M_1 \# M_2, \xi_1 \# \xi_2)$ by using the contact charts 
$\phi^1_1$ and $\phi^1_2$. 
Note that $k_1 \ge 2$ and $k_2 \ge 2$.
We can assume that $2 \le k_1 \le k_2$.
For $j \in \{ k_1+1, \dots, k_2 \}$ set $\cu^j = \cb^j_2(T)$.
For $i=1,2$ and $j \in \{ 2, \dots, k_1 \}$
let $\ck^j_i$ be the closure of $\cb^j_i(T)$.
For each $j \in \{ 2, \dots, k_1 \}$ the sets $\ck^j_1$ and $\ck^j_2$ are disjoint in $M_1 \# M_2$ 
and contained in $\cb^j_1$ resp.~$\cb^j_2$.
By Proposition~\ref{p:one} there exists a contact ball~$\cu^j$ in $(M_1 \# M_2,\xi_1 \# \xi_2)$ 
with $\ck^j_1 \cup \ck^j_2 \subset \cu^j$.
The sets $\cb^j_1(T), \cb^j_2(T)$ with $j \ge 2$ are then covered by the contact balls 
$\cu^2, \dots, \cu^{k_2}$.
Consider the ``finite neck'' 
$$
\cn (T) \,=\, \phi_1^1 (B_T \setminus B_{-1}) \cup \phi_2^1 (B_T \setminus B_{-1}).
$$
Since each $M_i$ is covered by $\cb_i^1(T), \dots, \cb_i^{k_i}(T)$ 
and since $\cb_1^j(T), \cb_2^j(T)$ with $j \ge 2$ are covered by $\cu^2, \dots, \cu^{k_2}$,
it will suffice to cover the set $\cn(T) \setminus \cu^2$ with one contact ball~$\cu^1$.
We distinguish two cases.

\m 
\ni
{\bf Case 1.} $2n+1 = 3$.

\s
\ni
Since $\cb^2_1(T)$ and $\cb^2_2(T)$ are not contained in $\cn (T)$, 
there exists an embedded smooth curve $C \subset \cu^2$ that starts in  
$\cb^2_1(T) \setminus \cn (T)$, ends in $\cb^2_2(T) \setminus \cn(T)$, and is such that
$C(T) := C \cap \cn (T)$ is connected.
Then $\cu^1 := \cn(T) \setminus C(T)$ is diffeomorphic to~$\RR^3$.
Since $(\cn, \xi_1 \# \xi_2)$ is contactomorphic to $(\RR^3 \setminus \{0\},\xi_{\rot})$ 
and hence to $\RR^3_{\st} \setminus \{0\}$, 
the set $(\cu^1, \xi_1 \# \xi_2)$ is therefore contactomorphic to a subset of $\RR^3_{\st}$ 
diffeomorphic to $\RR^3$. 
By Proposition~\ref{p:chart.3}, $\cu^1$ is a contact ball in $(M_1 \# M_2, \xi_1 \# \xi_2)$.
Since $\cu^1 \supset \cn(T) \setminus \cu^2$, the contact balls $\cu^1, \cu^2, \cu^3, \dots, \cu^{k_2}$ cover $M_1 \# M_2$.

\m 
\ni
{\bf Case 2.} $2n+1 \ge 5$.

\begin{figure}[ht] 
 \begin{center}
  \psfrag{R}{$T$}
  \psfrag{-R}{$-T$}
  \psfrag{2R}{$2T$}
  \psfrag{-2R}{$-2T$}
  \psfrag{cl}{$\cll$}
  \psfrag{g}{$\Gamma$}
  \psfrag{cn}{$\cn$}
  \leavevmode\epsfbox{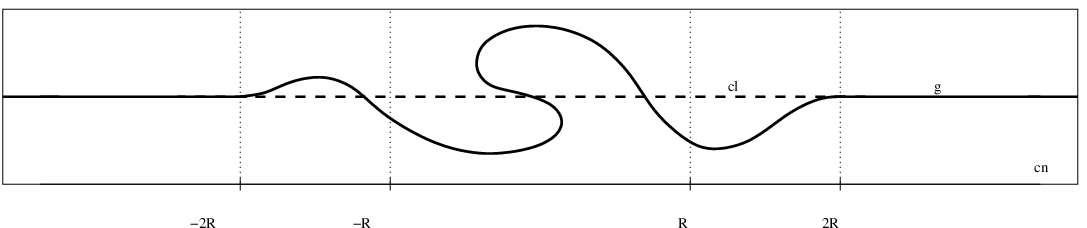}
 \end{center}
 \caption{The curves $\Gamma$ and $\cll$ in $\cn$.}
 \label{figure.gamma}
\end{figure}

\s
\ni
Consider the neck $\cn \subset M_1 \# M_2$ defined by~\eqref{def:neck},
and the lines $L = \left\{ (\00,\00,z) \in \RR^{2n+1} \mid z \ge 0 \right\}$ and $\cll = \phi(L) \subset \cn$.
Recall that $S_t = \pp B_t$ for $t \in \RR$.
We parametrise the lines $L$ and $\cll$ by $t \in \RR$ via $L(t) = L \cap S_t$ and 
$\cll (t) = \phi(L(t))$.
Then $\cll (t) \subset M_1 \setminus \cb_1^1(1)$ for $t \le -1$  
and  $\cll (t) \subset M_2 \setminus \cb_2^1(1)$ for $t \ge 1$. 
Since $\cb^2_1(T)$ and $\cb^2_2(T)$ are not contained in $\cn (T)$, 
there exists an embedded smooth curve $\Gamma \colon \RR \to \cn$ such that
\begin{equation}  \label{e:Gamma}
\left\{
 \begin{array}{l}
\Gamma \bigl( \, ]-T,T[ \, \bigr) \,\subset\, \cn(T) \cap \cu^2, \\[0.2em]
\Gamma \bigl( \, ]-2T,2T[ \, \setminus \, ]-T,T[ \, \bigr) 
                         \,\subset\, \cn(2T) \setminus \cn(T), \\[0.2em]
\Gamma(t) = \cll(t) \, \text{ for $|t| \ge 2T$},
\end{array}
\right.
\end{equation}
cf.~Figure~\ref{figure.gamma}.

%
%

\ni
Recall the contactomorphism 
$\phi \colon \left(\RR^{2n+1} \setminus \{0\}, \xi_{\rot} \right) \to (\cn, \xi_1 \# \xi_2)$.
We co-orient $\xi_1 \# \xi_2$ on~$\cn$ by the contact form $\ga := \phi_* \ga_{\rot}$.
A smooth curve $\gg \colon \RR \to \cn$ is {\it positively transverse}\/
if $\ga \left( \dot \gg(t) \right) = \gg^* \ga(t) >0$ for all $t \in \RR$.
The curve~$\cll$ is positively transverse and embedded.
Possibly after replacing $\Gamma$ with a $C^0$-close curve, 
we can assume that also~$\Gamma$ is positively transverse,
in view of a relative version of the $h$-principle 
(see \cite[Theorem~7.2.1]{EM}). 
After $C^\infty$ perturbing~$\Gamma$, if necessary, we can also assume that 
$\Gamma$ is embedded.
\begin{lemma}  \label{l:cpsi}
There exists a compactly supported contactomorphism~$\psi$ of $\cn$ such that $\psi(\Gamma) = \cll$.
\end{lemma}

\proof
Choose a smooth family $\Gamma_s \colon \RR \to \cn$, where $s \in [0,1]$, 
of smooth curves such that $\Gamma_0 = \Gamma$ and $\Gamma_1 = \cll$, and such that 
$\Gamma_s (t) \in \cn(2T)$ for all $s \in [0,1]$ and $|t| < 2R$.
Since $\Gamma_0$ and $\Gamma_1$ are positively transverse, we can apply a relative parametric 
$h$-principle
(see Theorem~7.2.1 and also Theorem~12.3.1 in \cite{EM}) 
and find a smooth family $\widetilde \Gamma_s \colon \RR \to \cn$ of 
positively transverse curves such that 
$\widetilde \Gamma_0 = \Gamma$ and $\widetilde \Gamma_1 = \cll$,
such that $\widetilde \Gamma_s (t) \in \cn(3T)$ for all $s \in [0,1]$ and $|t| < 3R$,
and such that $\widetilde \Gamma_s(t) = \cll(t)$ for $|t| \ge 3R$.
Each curve $\widetilde \Gamma_s$ is immersed.
Since $\dim M \ge 5$, we can perturb the family $\widetilde \Gamma_s$ 
such that each curve $\widetilde \Gamma_s$ is embedded.

We wish to extend the isotopy of curves $\widetilde \Gamma_s$ to a contact isotopy of~$\cn$ with support in~$\cn(4T)$.
Since each curve~$\widetilde \Gamma_s$ is a contact submanifold, 
the existence of such an isotopy can be easily deduced from the proof of the 
Isotopy Extension Theorem for contact submanifolds (Theorem~2.6.12 in~\cite{Gei-book}).
We prefer to give a direct argument, which is easier in our particular setting.
Define the time-dependent vector field $X_s$ along $\widetilde \Gamma_s$ by 
$$
X_s \left( \widetilde \Gamma_s(t) \right) \,=\, \frac{d}{ds} \widetilde \Gamma_s(t) .
$$
Fix $s \in [0,1]$. The Normal Form Theorem for transverse curves
(\cite[Example~2.5.16]{Gei-book}) asserts that there are
coordinates $(t,\xx,\yy)$ near $\widetilde \Gamma_s$ such that 
$$
\widetilde \Gamma_s (t) = (t,\00,\00) 
\qquad \text{and} \qquad 
\ga (t,\xx,\yy) = dt + \xx \,d\yy .
$$
In these coordinates, write
$$
X_s (t,\00,\00) \,=\, a_s(t) + \bb_s(t) \,\pp_{\xx} + \ccc_s(t) \,\pp_{\yy} .
$$
Define the smooth function $H_s$ near $\widetilde \Gamma_s$ by
$$
H_s (t,\xx,\yy) \,=\, a_s(t) + \ccc_s(t)\, \xx - \bb_s(t)\, \yy .
$$
The contact vector field $X_{H_s}$ of $H_s$ defined by~\eqref{e:contactHam}
equals $X_s$ along $\widetilde \Gamma_s$.
The coefficients $a_s(t), \bb_s(t), \ccc_s(t)$ vanish for $|t| \ge 3R$, 
and so $H_s(t,\xx,\yy) =0$ for $|t| \ge 3R$.
The coordinates $(t,\xx,\yy)$ can be chosen to depend smoothly on $s$.
After suitably cutting off the functions $H_s$, we therefore obtain a smooth function~$G_s$ on $\cn$ with support in $\cn(4T)$ 
such that $X_{G_s} = X_s$ along~$\widetilde \Gamma_s$.
The time~1~map~$\psi$ of the flow of $X_{G_s}$ is a contactomorphism 
of $\cn$ with support in~$\cn(4T)$ and such that 
$\psi \bigl( \Gamma (t) \bigr) = \cll (t)$ for all $t \in \RR$.
\proofend

\begin{lemma} \label{l:coverneck}
Any compact subset $\ck$ of 
$\left( \cn \setminus \cll, \xi_1 \# \xi_2 \right)$ 
can be covered by one contact ball in 
$\left( \cn \setminus \cll, \xi_1 \# \xi_2 \right)$.                          
\end{lemma}

\proof
The contactomorphism
$\phi \colon \left( \RR^{2n+1} \setminus \{0\}, \xi_{\rot} \right) \to \left( \cn, \xi_1 \# \xi_2 \right)$ 
from~\eqref{def:phi}
restricts to a contactomorphism
\begin{equation}  \label{e:LL}
\left( \RR^{2n+1} \setminus L, \xi_{\rot} \right) \to \left( \cn \setminus \cll, \xi_1 \# \xi_2 \right) 
\end{equation}
It therefore suffices to show that any compact subset $K$ of 
$\left( \RR^{2n+1} \setminus L, \xi_{\rot} \right)$ 
can be covered by one contact ball in  
$\left( \RR^{2n+1} \setminus L, \xi_{\rot} \right)$.
Choose an open neighbourhood $U$ of $K$ which is disjoint from~$L$
and such that
$$
(\xx, \yy,z) \in U \,\Longrightarrow\, (\xx,\yy,z') \in U 
\quad \text{ for all }\, z' \le z .
$$
Let $f \colon \RR^{2n+1} \to [0,1]$ be a smooth function such that $f |_L \equiv 0$ and $f |_U \equiv 1$.
The vector field $- \pp_z = (\00,\00,-1)$ preserves the contact form $\ga_{\rot}$ and has contact Hamiltonian $H \equiv -1$. 
The contact flow $\phi_{fH}^t = \phi_f^t$ preserves $\RR^{2n+1} \setminus L$,
and by the choice of $U$ we find $T>0$ such that $\phi_f^T(K) \subset \{ (\xx,\yy, z) \mid z<0\}$.
Choose $R$ so large that the open ball~$B$ of radius $R$ and centre~$(\00,\00,-R)$ covers $\phi_f^T(K)$.
In view of Proposition~\ref{p:contact.starshaped}, 
$B$ is a contact ball in $(\RR^{2n+1}\setminus L, \xi_{\rot})$.
Hence $\phi_f^{-T}(B)$ is a contact ball in $(\RR^{2n+1}\setminus L, \xi_{\rot})$ that covers $K$.
\proofend

Recall that we want to cover $\cn(T) \setminus \cu^2$ with one contact ball~$\cu^1$. 
In view of~\eqref{e:Gamma} we find a compact subset~$\ck$ of $\cn \setminus \Gamma$
such that $\cn(T) \setminus \cu^2 \subset \ck$.
With~$\psi$ the contactomorphism from Lemma~\ref{l:cpsi}, we then have $\psi(\ck) \subset \cn \setminus \cll$. 
By Lemma~\ref{l:coverneck}, there is a contact ball $\widetilde \cu^1$ in $\cn$ covering $\psi (\ck)$.
For the contact ball~$\cu^1 := \psi^{-1} \bigl( \widetilde \cu^1 \bigr)$ we therefore have 
$\cn(T) \setminus \cu^2 \subset \ck \subset \cu^1$. The contact balls $\cu^1, \cu^2, \dots, \cu^{k_2}$ cover $M_1 \# M_2$, and the proof of Theorem~\ref{t:connectedsum.d} 
is complete. 
\proofend

We conclude this section by proving 
\begin{proposition}  \label{p:Lcont}
$(\RR^{2n+1} \setminus L, \xi_{\rot})$ is contactomorphic to $\RR^{2n+1}_{\st}$.
\end{proposition}

Together with the contactomorphism~\eqref{e:LL} it follows that the ``neck without the line'' $\left( \cn \setminus \cll, \xi_1 \# \xi_2 \right)$ is contactomorphic to $\RR_{\st}^{2n+1}$. This strengthens Lemma~\ref{l:coverneck}.

\begin{figure}[ht] 
 \begin{center}
  \psfrag{xy}{$\xx, \yy$}
  \psfrag{z}{$z$}
  \psfrag{L}{$L$}
  \leavevmode\epsfbox{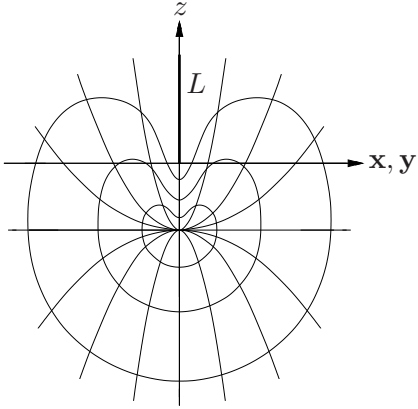}
 \end{center}
 \caption{The family $U_1 \subset U_2 \subset U_3 \dots$.} 
 \label{figure.hearts}
\end{figure}

\m
\ni
{\it Proof of Proposition~\ref{p:Lcont}.}
Since the contactomorphism~\eqref {e:cont} between 
$\left(\RR^{2n+1}, \xi_{\rot}\right)$ and $\RR^{2n+1}_{\st}$ preserves 
$L = \left\{ (\00,\00,z) \in \RR^{2n+1} \mid z \ge 0 \right\}$,
it suffices to show that $\left(\RR^{2n+1} \setminus L, \xi_{\st}\right)$ is contactomorphic to $\RR^{2n+1}_{\st}$.
Choose a sequence of bounded domains
$$
U_1 \,\subset\, U_2 \,\subset\, U_3 \,\subset\, \cdots
$$
such that for each $U_i$ the contact vector field $V (\xx,\yy,z) = (\xx,\yy, 2z+1)$ of $\RR^{2n+1}_{\st}$
is transverse to $\pp U_i$, 
such that $\overline U_i \subset U_{i+1}$,
and such that $\bigcup_{i=1}^\infty U_i = \RR^{2n+1} \setminus L$,
cf.~Figure~\ref{figure.hearts}.

%
%

\ni
By Proposition~\ref{p:contact.starshaped}, each set $(U_i, \xi_{\st})$ is contactomorphic to $\RR^{2n+1}_{\st}$.
Proposition~\ref{p:Lcont} is therefore a special case of the following
\begin{proposition}  \label{p:balls}
Let $U$ be a subset of $\RR^{2n+1}_{\st}$ that is the union $\bigcup_{i=1}^\infty U_i$ 
of bounded domains $U_i \subset \RR^{2n+1}_{\st}$ with the following properties:
$\overline U_i \subset U_{i+1}$,
and each $(U_i, \xi_{\st})$ is contactomorphic to $\RR^{2n+1}_{\st}$.
Then $(U, \xi_{\st})$ is contactomorphic to $\RR^{2n+1}_{\st}$.
\end{proposition}

\proof
We follow again Section~2.1 in~\cite{EG}.
Fix contactomorphisms $\gf_i \colon (U_i,\xi_{\st}) \to \RR_{\st}^{2n+1}$.
For $R>0$ let $B_R$ be the open ball in $\RR^{2n+1}$ of radius $R$ centred at the origin, 
and let $\overline{B_R}$ be its closure.
Set $\widetilde \gf_2 = \gf_2 \colon U_2 \to \RR^{2n+1}$.
Since $\overline U_1 \subset U_2$ is compact, we can find $R_1 \ge 1$ such that $\widetilde \gf_2 (U_1) \subset B_{R_1}$.
By the Contact Disc Theorem \cite[Theorem~2.6.7]{Gei-book},
applied to the contact embeddings  
$id$ and $\gf_3 \circ \widetilde \gf_2^{-1} \colon \overline {B_{R_1}} \to \RR_{\st}^{2n+1}$,
there exists a contactomorphism $\psi_1$ of $\RR_{\st}^{2n+1}$ such that
$\psi_1 \circ (\gf_3 \circ \widetilde \gf_2^{-1}) =  id$ on $\overline{B_{R_1}}$.
For the contactomorphism $\widetilde \gf_3 := \psi_1 \circ \gf_3 \colon U_3 \to \RR^{2n+1}$ we then have
$$
\widetilde \gf_3 = \widetilde \gf_2 \quad 
\text{ on }\, \widetilde \gf_2^{-1}(B_{R_1}) \supset U_1.
$$
Proceeding in this way we successively choose radii $R_i \ge i$ such that 
$\widetilde \gf_{i+1}(U_i) \subset B_{R_i}$ and construct contact embeddings 
$\widetilde \gf_{i+2} \colon U_{i+2} \to \RR^{2n+1}$ such that
$$
\widetilde \gf_{i+2} = \widetilde \gf_{i+1} \quad 
\text{ on }\, \widetilde \gf_{i+1}^{-1}(B_{R_i}) \supset U_i 
$$
for $i \ge 1$.
Since $U_1 \subset \widetilde \gf_2^{-1}(B_{R_1}) \subset U_2 \subset \widetilde \gf_3^{-1}(B_{R_2}) \subset \dots$
we have
$$
\widetilde \gf_{i+1}^{-1}(B_{R_i}) \subset \widetilde \gf_{i+2}^{-1}(B_{R_{i+1}})
\quad \text{for each }\, i
\qquad \text{and} \qquad
\bigcup_{i=1}^\infty \widetilde \gf_{i+1}^{-1}(B_{R_i}) = U .
$$
We can now consistently define a contact embedding $\Phi \colon U \to \RR^{2n+1}$
by 
$$
\Phi (u) := \widetilde \gf_{i+1}(u) \quad 
\text{ if }\, u \in \widetilde \gf_{i+1}^{-1} (B_{R_i}) \,\text{ for some }\, i \ge 1.
$$
Moreover, 
$$\Phi (U) \,\supset\, \widetilde \gf_{i+1} \left( \widetilde \gf_{i+1}^{-1} (B_{R_i}) \right) \,=\, B_{R_i} \subset B_i 
$$
for each $i \ge 1$, whence $\Phi (U) = \RR^{2n+1}$.
\proofend

\end{document}